\documentclass{article}

\usepackage[english]{babel}

\usepackage[letterpaper,top=2cm,bottom=2cm,left=3cm,right=3cm,marginparwidth=1.75cm]{geometry}
\usepackage[authoryear]{natbib}
\usepackage{booktabs}     
\usepackage{multirow}     
\usepackage{caption}      
\newtheorem{claim}{Claim}
\newtheorem{proposition}{Proposition}
\newtheorem{remark}{Remark}
\newtheorem{definition}{Definition}
\newtheorem{proof}{Proof}
\usepackage{algorithm}        
\usepackage{algpseudocode}    
\usepackage{amsmath,amssymb}  
\usepackage{subcaption}

\usepackage{amsmath}
\usepackage{graphicx}
\usepackage[colorlinks=true, allcolors=blue]{hyperref}

\title{Robust and Flexible Microtransit Design: Chance-Constrained Dial-a-Ride Problem with Soft Time Windows}
\author{
  Hongli Li\thanks{Lyles School of Civil and Construction Engineering, Purdue University, West Lafayette, IN, USA. Email: \texttt{li5125@purdue.edu}} \and
  Zengxiang Lei\footnotemark[1]\thanks{Email: \texttt{lei67@purdue.edu}} \and
  Xinwu Qian\thanks{Department of Civil and Environmental Engineering, Rice University, 77005. Email: \texttt{xq15@rice.edu}} \and
  Satish V. Ukkusuri\footnotemark[1]\thanks{Corresponding author. Email: \texttt{sukkusur@purdue.edu}}
}
\begin{document}
\maketitle

\begin{abstract}
Microtransit offers a promising blend of rideshare flexibility and public transit efficiency. In practice, it faces unanticipated but spatially aligned requests, i.e., passengers seeking to join ongoing schedules, leading to underutilized capacity and degraded service quality if not properly addressed. At the same time, microtransit needs to accommodate diverse passenger requests, ranging from routine errands to time-sensitive trips such as medical appointments. To meet these varying service expectations, incorporating time flexibility is essential. Nevertheless, existing models rarely account for both spontaneous and heterogeneous demand, limiting their effectiveness in real-world operations. We present a robust and flexible microtransit framework that simultaneously incorporates time flexibility and demand uncertainty through a formulation of a Chance-Constrained Dial-A-Ride Problem with Soft Time Windows (CCDARP-STW). We exploit the nonlinear chance constraints to model demand uncertainty with a controllable violation probability, ensuring high service reliability, and capture the time flexibility by soft time windows with penalized cost. A linear bounded-support relaxation through limited statistical information is developed to tackle the nonlinearity of the chance constraints. To efficiently solve the problem, we develop a tailored Branch-and-Cut-and-Price (BCP) algorithm enhanced with a probabilistic dominance rule. This rule improves computational performance by reducing the number of explored labels by 17.40\% and CPU time by 22.27\% in robust scenarios. A case study using real-world trip data from Chicago demonstrates the effectiveness of introducing time flexibility and robustness in microtransit. It achieves average savings of 11.55 minutes in travel time and 11.13 miles in distance compared to conventional microtransit. It also achieves the highest reliability (service rate of 96.46\%) among other robust approaches (i.e., normal distribution and ambiguity set) while maintaining computational tractability. 
\end{abstract}

\vspace{1em}

\section{Introduction}
Microtransit is a flexible and responsive transportation mode that combines the on-demand convenience of ridesharing with the operational efficiency of public transit. In low-density or underserved areas, it serves as a vital complement to existing transportation networks, particularly where conventional transit fails to meet community mobility needs\footnote{https://n-catt.org/guidebooks/microtransit-what-is-it-and-why-use-it-factsheet/}.

A key challenge in microtransit operations lies in the inherent uncertainty of travel demand. Empirical studies have shown that transportation demand forecasts are often inaccurate, with substantial discrepancies between predicted and realized demand~\citep{flyvbjerg2005accurate,peled2021quality}. In microtransit systems, this uncertainty can result in missed service opportunities or underutilized resources when systems rely solely on pre-booked or forecasted requests. For example, additional passengers seek to join ongoing shared trips from spatially concentrated locations such as hospitals, transit hubs, or commercial centers, leading to unanticipated demand at the planning stage but highly relevant operationally. \citet{wu2024prediction} demonstrate that incorporating spatial correlation in demand patterns in on-demand schedules can significantly improve service reliability. 

Another key challenge in microtrasnit operations stems from the diverse time sensitivity of passenger needs. As highlighted by \cite{ghimire2024policy}, an effective microtransit system must distinguish between urgent trips (e.g., for medical appointments or job interviews) and more flexible requests (e.g., shopping or routine errands). Modeling this heterogeneity through a mix of soft time windows and hard time windows allows for greater operational adaptability and improved alignment with user preferences. 

Addressing both demand uncertainty and time flexibility in a unified framework offers significant operational and service-level benefits in the microtransit system. Incorporating time flexibility not only meets user preferences but also enables the system to absorb more unanticipated requests due to its flexible schedules. In addition, accounting for demand uncertainty improves the system's ability to respond to different demand patterns and increase service reliability. Together, both features enhance the system's ability to serve a broader range of passengers, improve vehicle utilization, and maintain high service reliability under real-world conditions.   

Despite recent progress in developing adaptive microtransit models, the literature remains limited in integrating time flexibility and demand uncertainty within a single operational framework. Most existing studies focus on strategic planning or assume deterministic input data, overlooking the practical complexities of service delivery~\citep{miah2020barriers,rath2023microtransit}. While some works incorporate user heterogeneity or stochastic demand, they often treat time constraints as rigid or fail to model the spatial and temporal correlations of spontaneous requests~\citep{truden2021analysis,rossetti2023commuter}. Although \citet{cummings2024deviated} introduce a two-stage stochastic model for deviated fixed-route microtransit services, their approach does not capture variability in passenger time sensitivity. As such, there remains a critical need for models that jointly account for time flexibility and robustness to ensure feasible and efficient microtransit operations in uncertain environments.  

\begin{figure}[ht]
    \centering
    \includegraphics[width=1.0\linewidth]{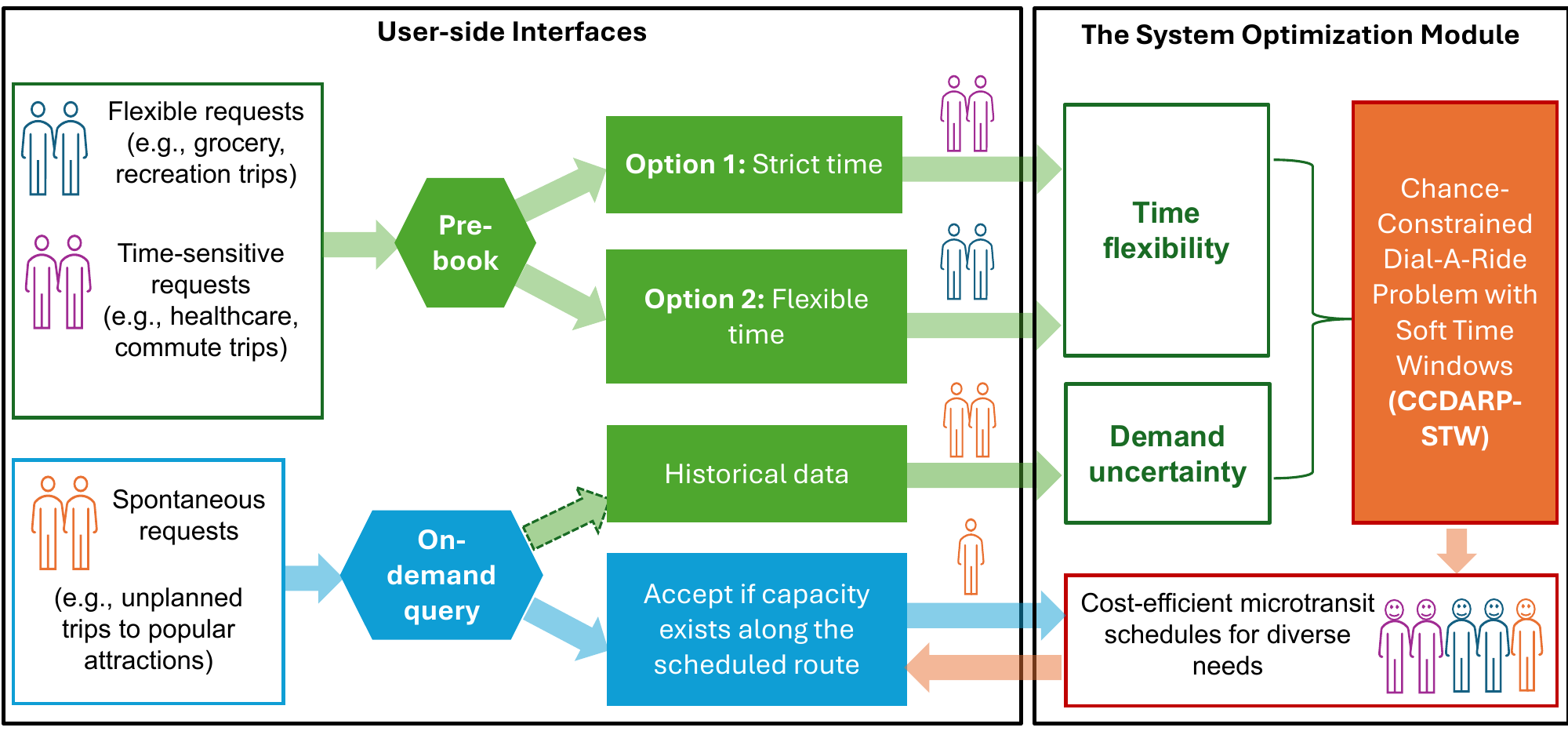}
    \caption{Robust and flexible microtransit framework}
    \label{fig:framework}
\end{figure}

This paper responds to the identified research gap by developing a pre-booked microtransit service framework that is both robust to demand uncertainty and adaptable to heterogeneous user needs. The system operates in multiple rounds and integrates both pre-booked and on-demand user interfaces, as illustrated in Figure~\ref{fig:framework}. For each service round, the underlying problem is formulated as a Chance-Constrained Dial-a-Ride Problem with Soft Time windows (CCDARP-STW). We introduce soft time windows to represent varying passenger urgency, enabling the system to delay flexible requests at a penalty while prioritizing critical trips and model the Demand uncertainty via chance constraints, ensuring service-level guarantees with a specified violation probability $\psi$ (e.g., 5\%). To enhance tractability, we reformulate the problem into a robust trip-based model and adopt a bounded-support linear relaxation based on Hoeffding's inequality, using only basic statistical information (mean and bounds).

To solve the model efficiently, we develop a tailored Branch-and-Cut-and-Price (BCP) algorithm. The pricing subproblem is addressed through a modified forward labeling algorithm that integrates soft time windows into label extension and pruning. To further improve scalability under uncertainty, we introduce a probabilistic dominance rule (PDR), which compares labels based on the probability of capacity violation. This rule enables early elimination of unpromising paths, reducing label exploration and computation time. Notably, the PDR is generalizable to other robust vehicle routing optimization problems involving probabilistic feasibility criteria. 

We evaluate our approach through comprehensive computational experiments on benchmark instances~\citep{Cordeau_2006} and a case study using real-world transportation data from Chicago. The BCP algorithm, augmented with the PDR, achieves a 17.40\% reduction in label exploration and a 22.27\% reduction in CPU time compared to standard dominance logic. In the case study, the proposed microtransit framework reduces vehicle usage by 67.67\% relative to ride-hailing and achieves savings of 11.55 minutes in time and 11.13 miles in distance compared to conventional microtransit. To guarantee a 95\% service level, it requires 26\% more vehicles than the deterministic scenario and achieves the highest service reliability, with an average service rate of 95.53\%, calculated as the number of successfully served requests divided by the total number of requests under demand uncertainty. Sensitivity analysis reveals that increasing time flexibility significantly reduces travel time and mileage up to a 50\% threshold, beyond which gains diminish. Among various demand modeling strategies, including normal distributions and ambiguity sets, the bounded-support approach offers the best trade-off, achieving the highest service reliability (96.46\%) while maintaining efficient resource usage. 

Section~\ref{sec:literature_review} reviews the relevant literature on microtransit systems, Dial-a-Ride Problems (DARP), and Branch-and-Cut-and-Price (BCP) algorithms. In Section~\ref{sec:problem_description}, we introduce the proposed modeling framework, including arc-based formulation and trip-based formulation. Section~\ref{sec:methodology} details the bounded-support approximation and the design of the tailored BCP algorithm. Section~\ref{sec:experiments} presents computational experiments that evaluate the performance and practical benefits of our proposed approach. Finally, Section~\ref{sec:conclusion} concludes the paper and outlines directions for future research.

\section{Literature Review}\label{sec:literature_review}
To fully appreciate the diverse contributions in modeling and solution approaches relevant to our problem, we organize the literature review into three parts. 

\paragraph{Operation of microtransit services} 
Microtransit has emerged as a promising mobility solution that combines the flexibility of ride-sharing with the efficiency of public transit. Positioned as a complement to a fixed-route system, micortransit is particularly well-suited for first-mile/last-mile connections, off-peak service coverage, and specialized transportation needs in low-density areas~\citep{lucken2019three,rossetti2023commuter}. When integrated effectively, microtransit can enhance accessibility, reduce travel time, and support broader urban sustainability goals by lowering dependence on single-occupancy vehicles~\citep{shaheen2016mobility}. Recent work by \citet{jacquillata2025microtransit} further investigates the design of microtransit services, highlighting the trade-offs and complementarities between fixed-line transit and on-demand mobility, and offering insights into when hybrid models are most effective.

A critical challenge in the deployment of microtransit services lies in managing uncertainty. Unlike fixed-route transit, microtransit operates in highly dynamic environments where passenger demand may deviate significantly from forecasts and evolve in real time. This makes robustness, a system's ability to maintain feasible and effective performance under variability, a key consideration. These advances in robust vehicle routing provide a valuable foundation for addressing uncertainty in microtransit operations. For instance, \citet{duan2021robust} propose a robust multiobjective optimization framework for the vehicle routing problem with time windows, which effectively balances cost efficiency and solution stability under uncertain travel conditions. \citet{de2019robust} further explore robustness in workforce-constrained vehicle routing, demonstrating how adaptable solutions can sustain service quality despite variability in operational parameters. Extending this line of work, \citet{vincent2023robust} develop a robust optimization model for vehicle routing with cross-docking, highlighting how accounting for demand uncertainty improves feasibility and operational resilience in logistics networks. While robust optimization has been widely applied to address uncertainty in logistics and vehicle routing problems, its integration into passenger-centered systems like microtransit remains limited~\citep{cummings2024deviated}. 

Yet robustness alone does not fully address the complexities of microtransit operations. Another remaining significant challenge in designing and operating microtransit systems is achieving time flexibility—the ability to provide reliable service to both time-sensitive and flexible passengers. The goal is to develop a microtransit system that can adapt to heterogeneous time requirements, offering precise scheduling for critical needs (e.g., healthcare appointments or meetings) while also accommodating passengers who tolerate moderate deviations. Building on the concept of heterogeneous travel patterns~\citep{ma2021user}, \citet{ghimire2024policy} emphasized the importance of distinguishing between these customer types to improve service responsiveness. \citet{JODIAWAN2025105173} further explored time flexibility in the flexible park-and-loop routing problem, showing how adaptive routing strategies enhance system performance across varied user expectations. 

Few studies have thoroughly addressed services that account for both time-sensitive and flexible customers and consider demand uncertainty. The integration of these factors into microtransit systems, with time flexibility and robustness, represents a novel challenge in creating more resilient and efficient transportation networks. This paper contributes to the development of a robust microtransit model that can handle both types of passenger needs while managing uncertainties in demand.

\paragraph{Dial a ride problem with time flexibility and demand uncertainty}
The Dial-a-Ride Problem (DARP) is a fundamental optimization model in demand-responsive transportation, where a fleet of vehicles must satisfy a set of pickup and drop-off requests while minimizing total travel cost and adhering to service constraints such as capacity, time windows, and ride time limits~\citep{cordeau2003tabu}. DARP models have been extensively used in paratransit, medical transport, and microtransit applications, where balancing efficiency and user satisfaction is essential.

In response to real-world variability and service heterogeneity, researchers have increasingly explored time flexibility through the use of soft time windows. Unlike hard time windows, which impose strict service periods, soft time windows permit slight deviations; early or late service incurs penalties rather than causing infeasibility~\citep{qureshi2009exact,slotboom2019impact}. This modeling approach is particularly important in systems that must accommodate both time-sensitive and flexible customers.

Several recent studies have advanced the use of soft time windows in vehicle routing and on-demand services. \citet{zhang2020multi} employed a multi-agent reinforcement learning approach to solve multi-vehicle routing problems with soft time windows, demonstrating that soft constraints enable more adaptive and efficient routing under uncertainty. Building on this line of work, \citet{zhang2023two} proposed a two-stage learning-based method for large-scale on-demand pickup and delivery services with soft time windows, further improving scalability and responsiveness in complex urban environments. In a more analytical context, \citet{chlebikova2023impact} examined how introducing soft ride time constraints affects the complexity of the Dial-a-Ride Problem (DARP), showing that while flexibility expands the solution space, it also increases computational challenges. In our work, we adopt a hybrid time-window structure, using soft time windows to serve flexible customers and hard time windows for those with strict timing needs, allowing tailored service and improved resource utilization.

Meanwhile, another crucial concern in DARP is uncertainty in demand. Traditionally, this has been handled through dynamic DARP, where routing decisions are made in real-time as requests arrive. \citet{coslovich2006two} proposed a two-phase insertion heuristic for dynamically arriving customers, and \citet{cordeau2007dial} provided a comprehensive review of dynamic DARP models and algorithms. These approaches emphasize reactive robustness, updating routes as uncertainty unfolds.

In contrast, our focus is on proactive robustness, scheduling under uncertainty in a static, pre-booked setting, where requests are not revealed dynamically but are uncertain at the scheduling stage. Here, we adopt a chance-constrained modeling approach, ensuring feasibility with high probability across stochastic demand realizations. The chance-constrained vehicle routing problem (CCVRP) has been widely studied for uncertain demands~\citep{ukkusuri2008linear}, including distributional robust models~\citep{sun2014distributionally,ghosal2020distributionally} and exact solutions using valid robust capacity inequalities~\citep{dinh2018exact}. Although a few studies have considered stochastic travel time in static DARP (e.g., \citet{ehmke2015ensuring}), the use of chance constraints for demand-side uncertainty in DARP remains limited.

To the best of our knowledge, no existing work integrates chance constraints over demand realizations into a static DARP with soft time windows. Our research addresses this gap by combining time flexibility and demand uncertainty in a unified framework, enabling more resilient pre-scheduled microtransit operations.

\paragraph{Branch-and-Cut-and-Price algorithm for solving complex DARP models}
As DARP models grow more complex, incorporating soft constraints, stochastic elements, and customer heterogeneity, traditional heuristics become insufficient. Metaheuristics such as tabu search~\citep{cordeau2003tabu}, genetic algorithms~\citep{jorgensen2007solving}, and granular search~\citep{KIRCHLER2013120} can solve large instances effectively but lack optimality guarantees and struggle with probabilistic feasibility. 

To solve advanced DARPs exactly, researchers have turned to decomposition-based methods such as Branch-and-Price (B\&P) and Branch-and-Cut (B\&C). The most powerful of these, the branch-and-cut-and-price (BCP) algorithm, combines the strengths of column generation, cutting planes, and branch-and-bound. It is particularly suited to set-partitioning or set-covering formulations, which naturally arise in DARP~\citep{Cordeau_2006}. The BCP framework iteratively generates promising routes (via pricing), strengthens relaxations (via cuts), and partitions the search space (via branching), enabling scalable exact solutions for large and constrained problems.

BCP has been successfully applied to various DARP extensions. \citet{luo2019two} tailored a two-phase BCP method for multi-depot DARPs with compatibility constraints. \citet{guo2025solving} implemented a BCP algorithm for a risk-aware DARP, optimizing routes while minimizing exposure risk. These studies illustrate how BCP can be adapted to domain-specific objectives and constraints.

Building on this foundation, we develop a tailored BCP algorithm to solve a DARP that incorporates both soft time windows and chance constraints. Our pricing problem is adapted to generate time-flexible routes while respecting probabilistic demand feasibility, and our master problem embeds both time-penalty logic and probabilistic coverage constraints. To the best of our knowledge, this is the first BCP framework that jointly handles time flexibility and demand uncertainty in a pre-scheduled microtransit setting. The algorithm allows efficient enumeration of robust and flexible service patterns while maintaining tractability for large-scale applications. 

\section{Problem Definition and Model Construction}\label{sec:problem_description}
To address the challenges of time flexibility and demand uncertainty in microtransit scheduling, we construct two mathematical models: an arc-based formulation and a trip-based robust formulation. The arc-based model explicitly represents vehicle routes over individual arcs and incorporates both soft and hard time windows, ride time constraints, and chance constraints on uncertain passenger loads. This formulation allows fine-grained control over temporal and load dynamics but can become computationally intensive for large instances. To overcome this, we develop an equivalent yet more tractable trip-based model, which shifts the decision focus from arcs to complete feasible trips. Leveraging the structural properties of DARP, we prove that the arc-level chance constraints can be equivalently aggregated at the trip level, allowing robust feasibility to be enforced across entire service patterns. This transformation enables a set-partitioning formulation that is more amenable to decomposition techniques like branch-and-price. Together, these models provide a flexible and scalable framework for designing resilient microtransit operations under demand uncertainty.
\subsection{Arc-based formulation}
The arc-based model $\mathcal{M}$ is formulated progressively. All notations used are summarized in Appendix \ref{sec:appendix_notation}. We consider a set of $n$ requests. Each request $i$ is associated with a pick-up node $i$ and a corresponding drop-off node $n+i$. Let $\mathcal{P}=\{1,2,\cdots,n\}$ denote the set of pick-up nodes and $\mathcal{D}=\{n+1,n+2,\cdots,2n\}$ the set of drop-off nodes. Each vehicle $k \in \mathcal{K}$ with capacity $M$, departs from the origin depot $0$ and ends its route at the destination depot labeled $2n+1$. The transportation network is modeled as a graph $\mathcal{G}=(\mathcal{N},\mathcal{A})$, where the node set is $\mathcal{N}=\mathcal{P}\cup \mathcal{D}\cup \{0,2n+1\}$, and the arc set is $\mathcal{A}=\{(i,j)~|~i,j \in \mathcal{N}\}$. Each node $i \in \mathcal{N}$ is associated with a time window $[e_i,l_i]$, where $e_i$ and $l_i$ denote the earliest and latest allowable service times, respectively. Requests are categorized into four types based on the flexibility of their time windows.
\begin{itemize}
    \item \textbf{Type (a)}: Fixed (hard) pick-up window and flexible (soft) drop-off window;
    \item \textbf{Type (b)}: Flexible pick-up window and fixed drop-off window;
    \item \textbf{Type (c)}: Fixed pick-up and drop-off windows;
    \item \textbf{Type (d)}: Flexible pick-up and drop-off windows.
\end{itemize}
Here, a \textit{flexible} (or soft) time window allows the latest service time $l_i$ to be extended if necessary, while a \textit{fixed} (or hard) time window must be strictly abided. Let $P \subseteq \mathcal{P}$ and $D \subseteq \mathcal{D}$ denote the subsets of pick-up and drop-off nodes, respectively, that are associated with flexible (soft) time windows. Define $N=P \cup D$ as the set of all flexible nodes. Accordingly, any node $i \in \mathcal{N} \setminus N$ is subject to a hard time window constraint that must be strictly satisfied. To enhance the competitiveness of microtransit, we incorporate \textit{load uncertainty} into the problem framework to accommodate unanticipated requests. Each request $i \in \mathcal{P}$ is associated with a random load $\tilde{w}_i$, satisfying the conservation condition $\tilde{w}_i=-\tilde{w}_{i+n}$, ensuring that each pick-up and its corresponding drop-off are load-balanced. The depot nodes have no associated load, i.e., $w_0=w_{2n+1}=0$. Furthermore, each request $i$ imposes a maximum allowable ride time $B_i$. The planning horizon is bounded by a global time limit $T$, requiring all vehicles to return to the destination depot before time $T$ after completing their assigned service.  

We define binary variables $x_{ij}^k$ to indicate whether vehicle $k \in \mathcal{K}$ traverses arc $(i,j) \in \mathcal{A}$, with $x_{ij}^k=1$ if the arc is used and 0 otherwise. Each arc $(i,j)$ incurs a travel cost $c_{ij}$. In addition to routing variables, we introduce several continuous side variables to capture the dynamics of the service process. Specifically, $L_i^k$ denotes the load carried by vehicle $k$ immediately after visiting node $i \in \mathcal{N}$, while $T_i^k$ represents the arrival time of vehicle $k$ at node $i$. The variable $D_i^k$ captures the ride time of request $i$ when transported by vehicle $k$. Furthermore, we introduce a delay cost function $H_i^k(\cdot)$ with a variable $y_i^k$ as the delay time experienced by vehicle $k$ when serving node $i \in N$. To address the possibility that some requests may remain unserved by the designated fleet, we incorporate a dummy vehicle $k=0$, which represents a third-party service provider. Assigning requests to this dummy vehicle incurs a significantly higher penalty cost, effectively discouraging such assignments except when no available vehicles can serve this request. 

Then, we introduce each set of constraints of model $\mathcal{M}$ step by step in the following paragraphs. The objective function~\eqref{model:objective function} minimizes the total operational cost, which includes both the travel cost incurred by vehicles and the penalties for delay.
\begin{equation}\label{model:objective function}
\mathcal{M}:=\min \quad \sum_{k \in \mathcal{K}} \sum_{i \in \mathcal{N}} \sum_{j \in \mathcal{N}} c_{ij} x_{ij}^k + \sum_{i \in \mathcal{N}} \sum_{k \in \mathcal{K}} H_i^k(y_i^k)  
\end{equation} 
The objective is minimized subject to a series of constraints that collectively ensure temporal feasibility, account for load uncertainty, enforce ride time limitations, and preserve network flow consistency.

We begin by modeling time-related constraints, including the representation of soft and hard time windows. Constraint~\eqref{constraint:time arc} tracks the cummulative travel time for each vehicle along the selected arcs. Constraint \eqref{constraint:flexible} captures both fixed and flexible time windows using an indicator function $\mathcal{X}_N(i)$, which activates the delay variable $y_i^k$ only if node $i$ belongs to the flexible set $N$. For depot nodes $0$ and $2n+1$, the time windows are fixed as $[0,T]$.
\begin{subequations}
    \begin{align}
        \label{constraint:time arc} T_j^k &\geq (T_i^k+s_i+t_{ij}) x_{ij}^k, \quad \forall i \in \mathcal{N}, j\in \mathcal{N}, k\in \mathcal{K} \\
         \label{constraint:flexible} e_i  &\leq T_i^k \leq l_i + \mathcal{X}_{N}(i) y_i^k, \quad  \forall i \in \mathcal{N}, k \in \mathcal{K} 
    \end{align}
\end{subequations}
To handle load uncertainty, we introduce a chance constraint in~\eqref{constraint:demand_chance}, which ensures that the vehicle load remains within capacity with high probability. Specifically, $1-\psi$ denotes the confidence level at which the vehicle's capacity is satisfied, while $\psi$ represents the acceptable risk of violation. This constraint guarantees that, with at least probability $1-\psi$, the load dynamics governed by the random variable $\tilde{w}_i$ do not exceed the vehicle's maximum capacity $M$. 
\begin{equation}\label{constraint:demand_chance}
    \text{Prob}\Bigg \{\begin{array}{rl}
         L_j^k & \geq (L_i^k + \tilde{w}_i)x_{ij}^k  \\
         L_j^k& \leq M  
    \end{array} \Bigg \} \geq 1- \psi, \quad \forall i,j \in \mathcal{N}, k \in \mathcal{K}
\end{equation}
Next, we enforce ride time limitations for each request. Constraint~\eqref{constraint:ride_time_1} calculates the ride duration of a request as the time difference between drop-off and pick-up. Constraint~\eqref{constraint:ride_time_2} ensures that this duration does not exceed the maximum allowable ride time $B_i$ for any request $i \in \mathcal{P}$. 
\begin{subequations}
    \begin{align} \label{constraint:ride_time_1}D_i^k&=T_{i+n}^k-T_i^k, \quad \forall i \in \mathcal{P} \\
    \label{constraint:ride_time_2}D_i^k &\leq B_i, \quad \forall i \in \mathcal{P}
    \end{align}
\end{subequations}
We then impose standard flow conservation and routing constraints to guarantee the feasibility of vehicle paths. Constraint~\eqref{constraint:flow conservation} ensures flow continuity at every passenger node. Constraints~\eqref{constraint:leaving once} and~\eqref{constraint:return once} enforce that each vehicle departs from the origin depot and returns to the destination depot exactly once. Constraint~\eqref{constraint:same vehicle} guarantees that the same vehicle visits both the pick-up and drop-off locations of a request. Finally, constraint~\eqref{constraint:single vehicle} ensures that each request is served exactly once across all vehicles.
\begin{subequations}
    \begin{align}
          \sum_{j \in \mathcal{N}}  x_{ji}^k -\sum_{j \in \mathcal{N}} x_{ij}^k=0, \quad &\forall i \in \mathcal{P}\cup \mathcal{D}, k \in \mathcal{K} \label{constraint:flow conservation}\\
    \sum_{j \in \mathcal{N}} x_{0j}^k=1, \quad &\forall k \in \mathcal{K}  \label{constraint:leaving once}\\
    \sum_{i \in \mathcal{N}} x_{i,2n+1}^k=1, \quad &\forall k \in \mathcal{K} \label{constraint:return once}\\
    \sum_{j \in \mathcal{N}} x_{ij}^k-\sum_{j\in N} x_{n+i,j}^k=0, \quad &\forall i \in \mathcal{P}, k \in \mathcal{K} \label{constraint:same vehicle}\\
    \sum_{k \in \mathcal{K}} \sum_{j \in \mathcal{N}} x_{ij}^k=1, \quad &\forall i \in \mathcal{P} \label{constraint:single vehicle}
    \end{align}
\end{subequations}
To complete the model, we impose standard variable domain constraints. Specifically, the assignment variables are binary, $x_{ij}^k \in \{0,1\}$m and all auxiliary variables, including $y_i^k, L_i^k, T_i^k$, and $D_i^k$, are non-negative for all $i,j \in \mathcal{N}$ and $k \in \mathcal{K}$. With these components, the arc-based formulation $\mathcal{M}$ is fully specified. To enhance computational efficiency, we adopt several reformulation techniques to make the arc-based model more compact.
\begin{remark}[ Implicit modeling of delay in soft time windows]
    Given the global time horizon $T$, the delay variable $y_i^k$ at each flexible node $i \in N$ is naturally bounded within the interval $[0,T-l_i]$. Leveraging this upper bound, we can omit the explicit modeling of $y_i^k$ and reformulate constraint~\eqref{constraint:flexible} as
$$e_i \leq T_i^k \leq l_i + \mathcal{X}_N(i)(T-l_i).$$
In this setting, we incorporate the delay penalties directly into a revised cost function by defining a new composite cost term $\hat{c}_{ij}^k$, which accounts for both travel cost and potential delay penalties at node $i$ and $j$. This composite cost is given by
$$
\hat{c}_{ij}^k=c_{ij}^k+H_i^k\Big((T_i^{k}-l_i)^+\Big) + H_i^k \Big ((T_j^{k}-l_j)^+\Big)
$$
Furthermore, we tighten the time window constraints for the depot nodes: for the origin depot, we impose $\min_{j \in \mathcal{P}} e_j \leq T_0^k \leq T$; and for the destination depot, we require $\min_{j \in \mathcal{D}} l_j \leq T_{2n+1}^k \leq T$ for all vehicles $k \in \mathcal{K}$.
\end{remark}

\begin{remark}[ Valid big-M reformulations for time and load consistency]
    We linearize the time constraint~\eqref{constraint:time arc} and the load constraints $L_j^k \geq (L_i^k+\tilde{w}_i)x_{ij}^k$ as below by introducing two large constraints $Q_1$ and $Q_2$. These constants act as valid big-M parameters when satisfying the following conditions~\citep{Cordeau_2006}: $Q_1 \geq \max_i\{0,l_i+\mathcal{X}_{N}(i)(T-l_i)+s_i+t_{ij}-e_j\},\ Q_2 \geq \min_i \{M, M+\tilde{w}_i\}.$
$$
\begin{array}{rl}
     T_j^k &\geq T_i^k + s_i + t_{ij} +Q_1(x_{ij}^k-1), \quad \forall i,j \in \mathcal{N}, k \in \mathcal{K} \\
     L_j^k &\geq L_i^k + \tilde{w}_i + Q_2(x_{ij}^k-1), \quad \forall i,j \in \mathcal{N}, k \in \mathcal{K}
\end{array}
$$
\end{remark}

Despite these reformulation efforts, solving the arc-based model remains computationally intensive due to the large number of decision variables and complexity introduced by time and uncertain load interactions across arcs. To address this issue, we introduce a trip-based formulation, which shifts the modeling perspective from individual arc traversals to feasible trips (i.e., sequences of requests served by a single vehicle). This alternative formulation significantly reduces the dimensionality of the model and enables a more tractable solution approach via column generation or decomposition methods.

\subsection{Trip-based formulation}
This section presents a trip-based robust model, denoted by $\mathcal{RM}$, which is derived from the arc-based formulation $\mathcal{M}$. The trip-based approach leads to a simpler model with fewer variables and constraints, making it easier to handle computationally. To formalize a trip-based model, we begin by defining feasible trips and then prove in Claim \eqref{claim:trip_based_chance_constraint} that the arc-based chance constraints can be equivalently written at the trip level. 

\begin{definition}\label{definition:route}
    A \textbf{trip} $r$, represented by a sequence of nodes  $\mathcal{I}_r=(i_0=0,i_1,\cdots,i_{m^r-1},i_{m^r}=2n+1)$, is said to be feasible in the microtransit scheduling context if it satisfies the following conditions: 
    \begin{itemize}
        \item a) For every request, the pick-up node $i \in \mathcal{P}$ appears in the sequence before its corresponding drop-off node $n+i \in \mathcal{D}$.
        \item  b) Each visited node $i$ satisfies the hard and soft time window constraint $e_i \leq T_i^k \leq l_i + \mathcal{X}_N(i) (T-l_i), \forall i \in \mathcal{N}$;
        \item c) Each request satisfies the maximum ride time constraint $D_i^k=T_{i+n}^k-T_i^k \leq B_i$. 
    \end{itemize}
\end{definition}

Definition \ref{definition:route} characterizes the structure of a feasible trip $r$ in the candidate trip set $\mathcal{R}$. A trip $r$ is deemed \textbf{robust} if it additionally satisfies the following trip-based chance constraint for each node $i \in \mathcal{I}_r$:
\begin{equation}\label{constraint:trip-based-formulation}
    \text{Prob}\Bigg \{\sum_{j \in O_i}\tilde{w}_j\leq M\Bigg \} \geq 1-\psi, \quad \forall i \in \mathcal{I}_r, r \in \mathcal{R}.
\end{equation}
where $O_i \subseteq \mathcal{P}$ denotes the set of requests that are on-board (i.e., picked up but not yet dropped off) at node $i$ along trip $r$. 

We claim that inequality \eqref{constraint:trip-based-formulation} provides an equivalent, aggregated, trip-level chance constraint. Although both formulations are equivalent in terms of conservativeness, they are different in modeling perspective. The arc-based chance constraint resolves violations on each arc, while the trip-based chance constraint considers the risk across the entire trip. The following claim formalizes this relationship and justifies the conservativeness of the trip-based formulation.
\begin{claim}\label{claim:trip_based_chance_constraint}
For a given trip $r$, arc-based formulation $\text{Prob}\{L_j+\tilde{w}_j\leq L_i \leq M, (j,i) \subset r\}$ for each $i \in I_r$ and trip-based formulation $\text{Prob}\{\sum_{j \in O_i}\tilde{w}_j \leq M\}$ for all $i \in I_r$ are equivalent, both satisfying being greater than $1-\psi$.  
\end{claim}
\begin{proof}
Let a trip $r$ be an ordered sequence of arcs from the origin depot to the destination depot. For example, suppose that $r$ consists of nodes and arcs such as:
$$
(0,2),(2,4),(4,1),(1,3),(3,5)
$$
where 0 and 5 are depots, 1 and 2 are pickup nodes, and 3 and 4 are the corresponding drop-off nodes. Suppose each request $j$ has a random demand $\tilde{w}_j \geq 0$, and let $L_i$ be the total load when the vehicle reaches node $i$. 

\paragraph{($\Rightarrow$) Arc-based implies trip-based:} Assume the arc-based constraints hold:
$$
\text{Prob} \{L_j+\tilde{w}_j \leq L_i \leq M\} \geq 1-\psi, \quad \forall (j,i) \subseteq r.
$$
We show that this implies:
$$
\text{Prob} \Bigg \{\sum_{j \in O_i} \tilde{w}_j \leq M\Bigg \} \geq 1-\psi, \quad \forall i \in I_r
$$
where  $O_i \subseteq \mathcal{P}$ denotes the set of requests that are on-board (i.e., picked up but not yet dropped off) at node $i$ along trip $r$ within a set of visited nodes $I_r$. By construction, the load at node $i$, denoted $L_i$, is $L_i=\sum_{j \in O_i \tilde{w}_j}$, i.e., the total uncertain load of customers still onboard just before node $i$. Since arc-based constraints are enforced along all arcs $(j,i)$ within trip $r$, they imply that the vehicle accumulates load consistently, the cumulative load does not exceed $M$, and the probabilistic condition is satisfied at each step. Hence, $L_i=\sum_{j \in O_i}\tilde{w}_j \leq M$ with probability at least $1-\psi$, for each $i \in I_r$. Therefore, the trip-based constraint holds.

\paragraph{($\Leftarrow$) Trip-based implies arc-based:} Now assume the trip-based constraint holds:
$$
\text{Prob} \Bigg \{\sum_{j \in O_i} \tilde{w}_j \leq M\Bigg \} \geq 1-\psi, \quad \forall i \in I_r.
$$
We want to show:
$$
\text{Prob} \{L_j+\tilde{w}_j \leq L_i \leq M\} \geq 1-\psi, \quad \forall (j,i) \subseteq r.
$$
Note that the arc-based constraint $L_j+\tilde{w}_j \leq L_i$ models the accumulation of load along the route. Since the load at node $i$ is by definition $L_i=\sum_{j \in O_i}\tilde{w}_j$, and since $L_j+\tilde{w}_j$ is a subset of that accumulation, this intermediate sum must be less than or equal to $L_i$.  Thus, if $\sum_{j \in O_i} \tilde{w}_j \leq M$ holds with probability $\geq 1- \psi$, then any subsum like $L_j+\tilde{w}_j \leq L_i \leq M$ also holds with at least that probability. Therefore, the arc-based constraints are satisfied under the trip-based formulation. 

\end{proof}
Claim \ref{claim:trip_based_chance_constraint} provides a theoretical foundation for translating arc-based chance constraints into an aggregated trip-based representation. This equivalence enables the formulation of a robust microtransit model that evaluates feasibility and uncertainty at the trip level rather than per arc, offering both computational and modeling advantages. Now, a robust trip-based (set-partitioning) model $\mathcal{TM}$ is displayed below
\begin{alignat}{2}
    \min \quad &  \sum_{r \in \mathcal{R}} C_r \lambda_r  \\
    \label{constraint:trip_assign}\mbox{s.t.} \quad &  \sum_{r \in \mathcal{R}} \beta_{i,r} \lambda_r = 1, \quad & i \in \mathcal{P}, \\
    \label{constraint:trip_single_vehicle}&\sum_{r \in \mathcal{R}}\lambda_r \leq |\mathcal{K}|, \\ 
    &\lambda_r \in \{0,1\}, \quad & r \in \mathcal{R},
\end{alignat}
where $C_r$ represents the total associated cost in a robust trip $r$ and $\beta_{i,r}$ is a binary parameter that request $i$ is served by the robust trip $r$. Constraint \eqref{constraint:trip_assign} represents that each request can be assigned to only one trip with one vehicle, and constraint \eqref{constraint:trip_single_vehicle} guarantees that the active trips should not exceed the maximum number of available vehicles. 

\section{Solution Methods}\label{sec:methodology}
To solve model $\mathcal{TM}$, we adopt a branch-and-cut-and-price (BCP) algorithm, which is recognized as a state-of-the-art method for solving this class of problems. A key challenge lies in constructing robust columns that satisfy the trip-based chance constraint \eqref{constraint:trip-based-formulation}. Additionally, generating valid and robust capacity cuts under uncertainty adds further complexity. Prior works such as \cite{dinh2018exact} and \cite{ghosal2020distributionally} have proposed methods to address trip-based chance constraints under assumptions like joint (independent or dependent) normal distribution and distributional ambiguity sets. However, in practice, the available information about uncertainty is often limited. To address this, we develop a bounded-support linear relaxation for the trip-based chance constraint, which leverages known properties such as the expectation, upper bound, and lower bound of the uncertain demand. This relaxation provides a tractable and conservative approximation, as detailed in Section~\ref{sec:bounded_support}. We also design a customized BCP algorithm that integrates time flexibility and demand uncertainty, as described in Section~\ref{sec:branch_and_cut_and_price}.

\subsection{Bounded support}\label{sec:bounded_support}
As mentioned before, one of the main challenges in solving the trip-based robust model $\mathcal{TM}$ is handling chance constraints under limited distributional information. In many real-world settings, we may only have partial knowledge about demand uncertainty, typically the expectation and bounded support-without access to the full distribution information. 

Previous studies have addressed this by assuming joint normal distributions or using distributional robust formulations that rely on known variance and expectation (e.g., \cite{dinh2018exact}, \cite{ghosal2020distributionally}). However, variance-based methods can be sensitive to distributional misspecification and may underestimate risk in practice. In contrast, we construct an uncertainty set based on the known upper and lower bounds of each request. This bounded support representation is not only distribution-free but also offers a more reliable and conservative approach to modeling robustness. It avoids reliance on variance and instead uses Hoeffding-type inequalities, which hold under both independence and weak dependence, as shown by \cite{janson2004large}. This makes it particularly well-suited for real-world applications with limited or noisy data. 

We assume that each request $i$ has an uncertain load $\tilde{w}_i$ with known expectation $\mu_i=\mathbb{E}[\tilde{w}_i]$ and bounded support $a_i \leq \tilde{w}_i \leq b_i$. Using Hoeffding's inequality, we obtain the following proposition, adapted from~\cite{shu2023humanitarian}, to conservatively approximate the trip-based chance constraint with a linear reformulation.
\begin{proposition}\label{proposition:demand_linear}
    The trip-based chance constraint \eqref{constraint:trip-based-formulation} can be conservatively approximated by the linear inequality \begin{equation}\label{eq:linear_demand_chance}
        \sum_{j \in O_i} \Big(\mu_j +\gamma (b_j-a_j)\Big) \leq M
    \end{equation}
    where $\gamma=\sqrt{\dfrac{\log (1/\psi)}{2}}$. 
\end{proposition}
\begin{proof}
See Appendix \ref{appendix:linear_constraint}
\end{proof}
The bounded-support linear approximation established in Proposition~\ref{proposition:demand_linear} not only simplifies the treatment of trip-based chance constraints but also enables the construction of robust valid inequalities for the master problem. In particular, this approximation provides a tractable means to evaluate whether a set of requests can be feasibly served within the vehicle capacity under uncertainty. This idea extends naturally to capacity cuts, which ensure that any subset of requests cannot be covered by fewer vehicles than necessary when demand uncertainty is considered. 

To define valid capacity cuts under uncertainty, let $q_r^{i,j}$ be the number of times edge $(i,j)$ appears in trip $r$, and let $\delta(S)$ denote the arc-cut set for subset $S \subseteq \mathcal{P}$, where:
$$
\delta(S)=\{(i,j) \in \mathcal{A}~|~i \in S, j \notin S, \ \text{or} \ i \notin S, j \in S\}.
$$
A robust capacity cut is given by 
$$
\sum_{r \in \mathcal{R}} \sum_{(i,j) \in \delta(S)} q_r^{i,j} \lambda_r \geq  2 \kappa_\psi, \quad \forall S \subseteq \mathcal{P},
$$
where:
$$\kappa_{\psi}:=\left \{\begin{array}{rl}
     1& \text{Prob}\Big\{\sum_{i \in S_p} \tilde{w}_i\leq M\Big \} \geq 1-\psi  \\
     2&  \text{otherwise}.
\end{array} \right.
$$
and $S_p$ denotes the set of requests served in subset $S$. By applying Proposition~\ref{proposition:demand_linear}, the chance constraint can be conservatively approximated by replacing $O_i$ with $S_p$ and applying the linear ineaulity:
$$
\sum_{i \in S_p} \Big ( \mu_i + \sqrt{\dfrac{1/\psi}{2}}\cdot (b_i-a_i)\Big) \leq M.
$$
These robust capacity cuts strengthen the master problem by explicitly enforcing feasibility under demand uncertainty. Rather than relying on deterministic thresholds, they incorporate probabilistic guarantees based on bounded-support information, without assuming full distributional knowledge. This ensures that the number of vehicles assigned to any subset of requests is sufficient even under worst-case realizations within known bounds. As a result, the cuts provide a practical and theoretically sound mechanism to enforce robustness in vehicle routing under uncertain and limited information, making them particularly valuable in real-world microtransit planning.

\subsection{Branch and cut and price algorithm}\label{sec:branch_and_cut_and_price}
Building on the bounded-support relaxation derived in Subsection~\ref{sec:bounded_support}, we now present a tailed BCP algorithm for the robust trip-based model $\mathcal{TM}$. We separate the linear relaxation of the model into a restricted master problem (RMP) within a partial set of columns and a pricing problem denoted as an elementary shortest path problem with resource constraints \citep{irnich2005shortest} due to exponential variables by the column generation procedures. In detail, we obtain the current optimal solution $\bar{\lambda}_r,r \in \bar{\mathcal{R}} \subset \mathcal{R}$ and the corresponding optimal dual solutions $\bar{p}_i,i \in \mathcal{P}$ and $\bar{q}$ by solving the restricted linear relaxation model. The pricing problem is to minimize the reduced cost $$z_r=C_r - \sum_{i \in \mathcal{P}} \beta_{i,r}\bar{p}_i-\bar{q}$$ to check whether there exist some robust columns with negative reduced costs.

\paragraph{Forward labeling algorithm} We devise a forward labeling algorithm to generate robust columns by considering time flexibility and demand uncertainty. We define a forward label $u$, representing a partial trip that starts from the origin depot $0$ and resides at the current node $n_u$, as below.
$$
u:=\Big(n_u, z_u, E_u,B_u, V_u, O_u,B_u^o, R_u^o(E_u),R_u^o(B_u^o), \Gamma_u \Big)
$$
where:
\begin{itemize}
    \item $n_u$ denotes the node where label $u$ resides,
    \item  $z_u$ is the reduced cost after serving node $n_u$,
    \item $E_u$ is the earliest service time and $B_u$ is the latest service time,
    \item $V_u$ is the set of completed requests along the partial trip until node $n_u$,
    \item $O_u$ is the set of open requests, where request $a \in O_u$ has been picked up but not dropped off,
    \item $R_u$ is the visited route until the resided node $n_u$,
    \item $B_u^o$ represents a feasible start-of-service time at node $n_u$ for an open request $a \in A_u$,
    \item $R_u^o(\tau)$ denotes the possible dropoff time for request $o \in O_u$ as a function of start-of-service time $\tau$. 
\end{itemize}
If $E_u \geq B_u^o$, it implies that $R_u^o(E_u)$ is bounded by the latest pickup time window or the dropoff time window $l_{n_u}$.

The violation probability $\Gamma_u$ is calculated using Proposition \ref{proposition:demand_linear} as:
\begin{equation}\label{equa:gamma_u}
\Gamma_u=\text{exp}\Bigg(-2\Big(\dfrac{M-\sum_{j \in O_u}\mu_j}{\sum_{j \in O_u}(b_j-a_j)}\Big)^2\Bigg)    
\end{equation}
when $O_u \neq \emptyset$ and $\Gamma_u=0$ when $O_u=\emptyset$.

To construct a feasible trip, we use the following resource extension functions to extend the label $u$ to a new label $u'$. Regarding time flexibility, we distinguish two cases in the next node $n_{u'}$, let $E_{u'}=E_u+s_{n_u}+t_{n_u n_{u'}}$,
\begin{itemize}
    \item If the node is critical and $E_{u'}>l_{n_{u'}}$, we reject this extension.
    \item If the node is flexible and $E_{u'}>l_{n_{u'}}$, we accept the extension with an added penalty cost  $c_p=d\cdot (E_{u'-}l_{n_{u'}})$.
\end{itemize}
$\Gamma_{u'}$ is updated using the same formula \eqref{equa:gamma_u} with $O_{u'}$ replacing $O_u$, and we ensure $\Gamma_{u'} \leq \psi$. Other updates-including the completed request set $V_{u'}$, the open request set $O_{u'}$, ride time information $B_{u'}^o$, $R_{u'}^a(E_{u'})$, and $R_{u'}^o(B_{u'}^o)$ following existing works such as~\citet{gschwind2018bidirectional},\citet{guo2025solving}, and \citet{yang2023exact}, until the label reaches the destination depot and forms a robust trip.

To accelerate the labeling process and ensure optimality, we introduce the following strong dominance rule.
\begin{proposition}[Probabilistic dominance rule]\label{proposition:strong_dominance}
    Given two robust labels $u_1$ and $u_2$ residing at the same ending node $n_1=n_2 \in \mathcal{N}$, label $u_1$ dominates $u_2$ if the following conditions are satisfied: 1) $z_1 \leq z_2$, $E_1 \leq E_2$, $V_1 \subseteq V_2$, and $O_1 \subseteq O_2$; 2) $R_1^o(E_1)-E_1 \geq R_1^o(E_2)-E_2$ for $o\in O_1$ and $R_1^o(B_1^o)\geq R_2^o(B_2^o)$ for $o \in O_1$; 3) $\Gamma_1 \leq \Gamma_2$. 
\end{proposition}
\begin{proof}
See Appendix \ref{appendix:proof_dominance}
\end{proof}
\begin{remark}[Generalization of probabilistic dominance rule]
    The proposed dominance rule incorporates a probabilistic measure, specifically, the violation probability $\Gamma$, into the label comparison procedure. This constitutes a probabilistic dominance rule that prioritizes not only cost and temporal feasibility but also reliability under uncertainty. In our case, $\Gamma$ quantifies the upper bound on the risk of violating capacity constraints due to uncertain demands, derived using Hoeffding’s inequality. By comparing labels based on this violation probability and favoring those with lower $\Gamma$, the algorithm achieves a probabilistically conservative search strategy. It provides a controlled price of robustness: the algorithm may retain labels with slightly worse reduced cost if they provide stronger guarantees of feasibility. This trade-off is crucial in practical settings where service reliability is prioritized over marginal gains in optimality.

Importantly, the probabilistic dominance rule is not limited to the context of robust microtransit or Dial-a-Ride Problems. It is generalizable to other labeling-based algorithms for stochastic or robust optimization problems, such as: Robust vehicle routing and inventory routing problems, where probabilistic constraints on load, time windows, or service levels are present; Stochastic shortest path or resource-constrained path problems, where uncertain travel times or demands impact feasibility; Multi-commodity network flow problems under uncertainty, where chance constraints on arc capacities or delivery reliability are imposed.

To apply this rule in a general setting, we offer a generalized template
\begin{enumerate}
    \item Define the relevant probabilistic constraint (e.g., reliability, feasibility, coverage) based on domain requirements;
    \item Quantify the violation risk using distributional knowledge or concentration inequalities such as Markov, Chebyshev, or Hoeffding bounds;
    \item Embed this risk metric into the label definition and construct dominance rules that compare both deterministic resources and probabilistic robustness.
\end{enumerate}
\end{remark}
This strategy enables the integration of uncertainty into the dominance logic, thereby guiding label-based algorithms toward feasible and resilient solutions while maintaining computational tractability. It offers a principled framework to unify robust decision-making and efficient search in complex, uncertain environments.

\paragraph{Branching and cutting} To strengthen the solution process within the BCP framework, we implement two types of branching strategies alongside a set of classical and robust cutting planes.
\begin{itemize}
    \item The first branching rule targets the number of vehicles in use. When the LP relaxation yields a fractional number of active routes, we define the fractional total as $\bar{X}v = \sum{r \in R'} \bar{X}r$, where $R'$ denotes the set of selected routes. We then generate two subproblems: one enforces an upper bound on the total number of selected routes, $\sum{r \in R'} X_r \leq \lfloor \bar{X}v \rfloor$, and the other enforces a corresponding lower bound, $\sum{r \in R'} X_r \geq \lceil \bar{X}_v \rceil$. After branching, the reduced cost of each route is adjusted as $z_r \leftarrow z_r - \alpha$, where $\alpha$ is the dual variable associated with the new constraint. This branching strategy is particularly effective in accelerating convergence when the LP solution is fractional at the route level.
    \item The second branching rule operates on arc flows between node partitions. Specifically, we leverage the structure of the arc-based relaxation by identifying a cut set $\delta(O)$ defined as the set of arcs $(i, j)$ such that $i \in O$ and $j \in \mathcal{N} \setminus O$, where $\mathcal{N}$ is the full node set. We select a subset $O$ such that the sum of arc variables across the cut, $\sum_{(i,j) \in \delta(O)} \bar{x}{ij}$, is approximately 1.5, indicating a highly fractional flow. Two branches are then created: one restricts the outflow across the cut to at most 1, and the other requires it to be at least 2. These are implemented as route-level constraints by aggregating over the associated route indicators $\beta{ij,r}$, resulting in the constraints $\sum_{(i,j) \in \delta(O)} \beta_{ij,r} X_r \leq 1$ and $\sum_{(i,j) \in \delta(O)} \beta_{ij,r} X_r \geq 2$ for all relevant cut sets $O \in \Omega$. The reduced cost is correspondingly updated as $z_r' \leftarrow \sum_{O \in \Omega} q_O \sum_{(i,j) \in \delta(O)} \beta_{ij,r}$, where $q_O$ is the dual variable associated with the cut constraint.
\end{itemize}
To further tighten the master problem and eliminate infeasible or inefficient routes, we incorporate several classes of valid inequalities. In addition to robust capacity cuts discussed in Subsection~\ref{sec:bounded_support}, we apply infeasible path inequalities, adapted from \citet{Cordeau_2006}, to remove routes that violate precedence or time window constraints. These cuts are particularly useful in enforcing temporal feasibility when the route structure is otherwise admissible. Furthermore, we employ two-path cuts as described in \citet{ropke2009branch}, which are used to exclude combinations of trips that, when scheduled together, lead to infeasible solutions. These cutting planes work collectively with the branching rules to ensure integrality, improve convergence speed, and maintain feasibility under both deterministic and uncertain conditions.

\paragraph{Algorithm framework} After developing the forward labeling algorithm to solve the pricing problem under time flexibility and demand uncertainty, and establishing both strong dominance rules and branching strategies to guide the search, we now outline the complete BCP framework in Algorithm \ref{alg:bcp}. The algorithm iteratively solves a restricted master problem (RMP) using a subset of feasible and robust trips. At each iteration, dual prices are extracted to guide the pricing subproblem, which identifies new promising routes with negative reduced cost using the forward labeling procedure. If such columns are found, they are added to the RMP. If no improving columns remain, the algorithm applies valid cutting planes, such as infeasible path and two-path inequalities, to further tighten the formulation. When neither new columns nor effective cuts are available, branching is performed to enforce integrality. This integrated framework balances column generation, constraint separation, and branching within a unified loop, allowing us to efficiently solve large-scale, robust, and flexible microtransit problems to optimality.
\begin{algorithm}[H]
\caption{Branch-and-Cut-and-Price Algorithm for \( \mathcal{TM} \)}
\label{alg:bcp}
\begin{algorithmic}[1]
    \State Initialize RMP with feasible columns \( \bar{\mathcal{R}} \subset \mathcal{R} \)
    \Repeat
        \State Solve LP relaxation of RMP to get \( \bar{\lambda}, \bar{p}, \bar{q} \)
        \State Solve pricing subproblem via forward labeling to find new columns
        \If{new columns found with \( z_r < 0 \)}
            \State Add new columns to RMP
        \Else
            \State Apply cutting planes (e.g., infeasible path, 2-path cuts)
            \If{no more valid cuts found}
                \State Apply branching (on vehicle count or arc flow)
            \EndIf
        \EndIf
    \Until{all columns and constraints explored and integer solutions found}
\end{algorithmic}
\end{algorithm}

\section{Computational Experiments and Case Study}\label{sec:experiments}
This section presents computational results to evaluate both the algorithmic behavior and the practical benefits of our proposed approaches. The experiments are divided into two parts:
\begin{itemize}
    \item Section~\ref{sec:algorithm_comparison} focuses on algorithmic comparisons using benchmark instances from~\cite{Cordeau_2006}. We assess the impact of moving from the classical DARP to more advanced formulations (flexible and robust DARP) on solution quality and computational efficiency. We also evaluate the influence of our proposed probabilistic dominance rule versus traditional dominance logic to understand how changes in labeling strategy affect the performance of the BCP algorithm in solving complex, uncertainty-aware models. Furthermore, we quantify the computational and operational costs associated with introducing time flexibility and robustness. Specifically, we examine how increasing time flexibility expands the labeling space, and how lowering the violation probability $\psi$ imposes higher resource demands in terms of vehicle capacity and fleet size. These analyses reveal important trade-offs and show that our BCP algorithm, supported by effective dominance rules and conservative scheduling, maintains strong performance across increasingly realistic and complex DARP settings.
    \item Section~\ref{sec:case_study} presents a case study based on real-world ride-hailing trip data from Chicago (2023-2024), focusing on a representative urban community. In this section, we compare five service strategies, ride-hailing, traditional microtransit, flexible microtransit, robust microtransit, and flexible and robust microtransit, to assess their practical implications for service efficiency, reliability, and resource usage. Beyond the direct comparison of service models, we conduct a series of sensitivity analyses to evaluate the impact of key design features. Specifically, we analyze how varying levels of time flexibility, violation probability $\psi$, and different demand modeling assumptions affect system performance. These analyses aim to quantify the trade-offs and operational benefits introduced by our proposed enhancements, offering guidance for deploying microtransit solutions under real-world uncertainty and complexity.  
\end{itemize}

\subsection{Algorithm Performance Comparison}\label{sec:algorithm_comparison}
We evaluate our algorithm using the standard benchmark datasets for DARP\footnote{\url{http://neumann.hec.ca/chairedistributique/data/darp/branch-and-cut/}}. We select the second set of instances, which are more appropriate for microtransit scenarios where each request may involve multiple passengers. In contrast, the first set of instances only demands one passenger for each request. Vehicles are heterogeneous with a fixed capacity of 6, and each request is constrained by a maximum ride time of 45 units. For example, the instance ``b2-16'' represents 2 vehicles serving 16 requests. The following columns are shared across all models in Table ~\ref{tab:time_sensitivity} and Table~\ref{tab:dominance_rule}:
\begin{itemize}
    \item ``Instance'': Name of the benchmark file;
    \item ``Bound'':  Optimal objective value for the corresponding model (equal to total travel time);
    \item ``CPU'': Total computation time in seconds.
    \item ``$N_p$'': Number of labels (i.e., feasible routes) explored during the pricing process in the BCP algorithm.
\end{itemize}

\paragraph{Algorithmic performance comparison on time flexibility and robustness} Table~\ref{tab:time_sensitivity} presents a comparative analysis across three models: classical DARP (C-DARP) with vehicle capacity 6, time-flexible DARP (TF-DARP) with the same capacity, and robust DARP (R-DARP) with increased capacity of 13. The C-DARP model requires all requests to follow strict hard time windows. The TF-DARP model introduces flexibility by allowing a portion of service nodes to deviate from their designated time windows at a cost.

Specifically, TF-DARP applies a flexibility ratio of 50\%, meaning half of the total service nodes (i.e., pick-up and drop-off nodes) are allowed to violate their time constraints. For instance, in a 16-request instance with 32 service nodes, 16 nodes are assigned soft time windows while the other 16 follow hard constraints. To account for time deviation, we adopt a linear delay cost function defined as $H_i^k=d_i^k \times y_i^k$, where $d_i^k=1$ for each request $i$ assigned to vehicle $k$. This formulation captures a balanced trade-off between flexibility and service punctuality. For alternative penalty structures and parameter settings tailored to specific instances, we refer readers to \citet{qureshi2009exact,slotboom2019impact}. It is worth noting that the choice of penalty cost function does not affect the formulation of our model, the design of our algorithm, or the resulting performance, as it remains an internal component of the optimization process. 

In the R-DARP model, we incorporate demand uncertainty by treating the pickup load at each node as a random variable with bounded support: the nominal load equals the mean demand, and its uncertainty interval spans from mean$-1$ to mean$+1$. A violation probability of $\psi=0.01$ is imposed to ensure high service reliability. To maintain feasibility under this conservative setting, we increase the vehicle capacity from 6 to 1 without changing the number of available vehicles.
\begin{table}[ht!p]
    \centering
    \begin{tabular}{c|ccc|ccc|ccc}
    \toprule
         & \multicolumn{3}{c}{C-DARP (6)} & \multicolumn{3}{|c}{TF-DARP (6)} & \multicolumn{3}{|c}{R-DARP (13)}  \\
    \hline
         Instance & Bound & CPU (s) & $N_p$ & Bound & CPU (s) & $N_p$  & Bound & CPU (s) & $N_p$ \\
    \hline
         b2-16 & 309.41 & 0.67 &392  & 309.41 & 1.21 & 436  & 309.41 & 3.29 & 342  \\
         b2-20 & 332.64 & 0.76 & 519  & 332.64 & 1.46 & 633  & 332.64 & 5.46 & 581   \\
         b2-24 & 444.71 & 39.26 & 922  & 444.71 & 44.91 & 1015  & 444.71 & 219.62 & 827  \\
         b3-18 & 301.64 & 1.32 & 248 & 301.64 & 1.19 & 289  & 301.64 & 3.98 & 218 \\
         b3-30 & 527.41 & 14.42 & 1206 & 527.41 & 37.6 & 2178  & 527.41 & 101.32 & 1237 \\
         b3-36 & 603.79 & 99.96 & 2255  & 603.79 & 220.02 & 4127  & 603.79 & 538.45 & 2482  \\
         b4-16 & 296.96 & 0.16 & 123  & 296.96 & 0.17 & 132  & 296.96 & 0.38 & 106  \\
         b4-24 & 371.41 & 2.07 & 562  & 371.41 & 5.31 & 983  & 373.56 & 36.93 & 586   \\
         b4-32 & 494.82 & 98.50 & 1319 & 494.82 & 235.25 & 2637  & 494.82 & 464.79 &2139 \\
         b4-40 & 655.36 & 127.76 & 1659  & 655.36 & 265.43 & 3428  & 655.52 & 574.86 & 2598  \\
         b5-40 & 613.72 & 113.10 & 1314 & 613.72 & 169.68 & 1389  & - & - & -  \\
         b6-48  &713.79 & 335.28 & 1318   & 713.79 & 342.91 & 1327   & - & - & - \\
         b6-60 & 858.88 & 747.59 & 2656   & 858.88 & 869.25 & 2680   & - & - & - \\
         \hline
         Avg. & - & 121.60 & 1115 & - & 168.80 & 1635 & - & 194.91 & 1112 \\
    \bottomrule
    \end{tabular}
    \caption{Algorithm Performance Comparison on Time Flexibility and Robustness}
    \label{tab:time_sensitivity}
\end{table}
The results in Table~\ref{tab:time_sensitivity} highlight several key findings:
\begin{enumerate}
    \item \textbf{C-DARP optimality}: Our BCP algorithm successfully solves all C-DARP instances to optimality, as indicated in the ``Bound'' column. Notably, TF-DARP yields the same objective values across all instances. This consistency suggests that, under current resource constraints (vehicle capacity and fleet size), C-DARP solutions remain optimal even when time flexibility is introduced. However, TF-DARP explores more feasible routes (column $N_p$) by allowing controlled time deviations, resulting in higher CPU times.
    \item \textbf{R-DARP feasibility and slight deviations}: With expanded vehicle capacity, R-DARP achieves the same optimal schedules as C-DARP in most cases. Minor deviations are observed in instances ``b4-24'' and ``b4-40'', where the objective values slightly increase. This can be attributed to the probabilistic dominance rule: the use of a violation threshold $\psi$ accelerates computation but may lead to marginal suboptimality in rare cases.
    \item \textbf{Scalability and performance}: The CPU times (column ``CPU (s)'') confirm the algorithm's efficiency even for more complex models like TF-DARP and R-DARP, which incorporate soft time windows and chance constraints. The growth in CPU time and number of labels remains reasonable and scalable.
    \item \textbf{Overall comparison}: The average statistics in the bottom row reinforce the scalability and robustness of our BCP algorithm. Although TF-DARP and R-DARP require more computational resources, they remain tractable and effective, validating the practical applicability of our framework for flexible and robust microtransit scheduling.
\end{enumerate}

\paragraph{Algorithmic performance comparison on dominance rules} To see the impact of our proposed probabilistic dominance rule (PDR) on solving robust models, we compare two variant of a time-flexible and robust DARP (TFR-DARP): one using a standard dominance rule (SDR), and the other applying the PDR. The TFR-DARP model integrates both time flexibility and demand uncertainty, as introduced separately in TR-DARP and R-DARP.

The primary distinction between SDR and PDR lies in how they handle the label feature $\Gamma$, an estimated violation probability. Under standard dominance, all labels satisfying $\Gamma \leq \psi$ are accepted. In contrast, the probabilistic dominance rule further compares labels by their respective violation probabilities: a label with a smaller $\Gamma$ is considered more promising and dominates another, allowing us to prune less robust labels early. This prioritization encourages more conservative routing paths under uncertainty and improves computational efficiency.
\begin{table}[h]
    \centering
    \begin{tabular}{c|ccc|cccc}
    \toprule
         & \multicolumn{3}{|c}{TFR-DARP without PDR} & \multicolumn{4}{|c}{TFR-DARP with PDR}  \\
    \hline
        & Bound & CPU (s) & $N_p$ & Bound & CPU (s) & $N_p$ & LR (\%) \\
    \hline
         b2-16  & 309.41 & 2.43 & 425 & 309.41 & 2.39 & 392 & 7.76 \\
         b2-20 & 332.64 & 4.58 & 646 & 332.64 & 4.34 & 607 & 6.04 \\
         b2-24  & 444.71 & 218.89 & 903 & 444.71 & 208.69 & 852 & 5.65  \\
         b3-18 & 301.64 & 3.51 & 255 & 301.64 & 3.27 & 211 & 17.25\\
         b3-30 & 527.41 & 183.18 & 2339 & 527.41 & 97.31 & 1704 &27.14 \\
         b3-36 & 603.79 & 885.96 & 3592 & 603.79 & 529.41 & 2979 & 17.01\\
         b4-16 & 296.96 & 0.37 & 126 & 296.96 & 0.28 & 112 & 11.11\\
         b4-24 & 371.41 & 36.73 & 604 & 373.56 & 34.83 & 567 & 6.13\\
         b4-32  & 494.82 & 708.76 & 3467 & 494.82 & 654.99  & 2939 & 15.23\\
         b4-40 & 655.36 & 793.42& 3967 & 655.52 & 670.16 & 3112 & 21.55\\
         \hline
         Avg. & - & 283.78 & 1632 & - & \textbf{220.57} & \textbf{1348} & \textbf{17.40}\\
    \bottomrule
    \end{tabular}
    \caption{Algorithm Performance Comparison on Probabilistic Dominance Rule}
    \label{tab:dominance_rule}
\end{table}
Table~\ref{tab:dominance_rule} summarizes the computational results of solving TFR-DARP with and without the probabilistic dominance rule. Column LR (\%) captures the reduction rate in label exploration achieved by PDR. It is computed as:
$$
\text{LR}=\dfrac{N_{p}^{\text{SDR}}-N_{p}^{\text{PDR}}}{N_{p}^{\text{SDR}}}\times 100\%
$$
where $N_p^{\text{SDR}}$ and $N_p^{\text{PDR}}$ represent the number of labels explored without and with PDR, respectively.

The results demonstrate that applying PDR significantly reduces the number of explored labels across all instances. On average, PDR achieves a 17.40\% reduction in label generation and a 22.27\% CPU time saving compared to using SDR. This efficiency gain translates into reduced computation time as well, without compromising optimality in most instances. While slight increases in the objective value are observed in ``b4-24'' and b4-40'', this trade-off is expected due to the more aggressive pruning of labels under PDR, which favors conservative and computationally efficient paths over exhaustive enumeration.

In summary, PDR provides a principled and effective mechanism for improving the tractability of solving robust DARP variants by filtering out less promising paths based on their violation probabilities, significantly accelerating the BCP algorithm with minimal impact on solution quality.

\paragraph{Price of time flexibility and robustness} We assess the computational and operational costs associated with introducing time flexibility (via soft time windows) and robustness (via chance constraints under demand uncertainty) in the DARP model solved by the proposed BCP algorithm. These two features enhance service adaptability and reliability, but also increase the computational load and resource requirements.
\begin{figure}[thbp]
    \centering
    \includegraphics[width=0.6\linewidth]{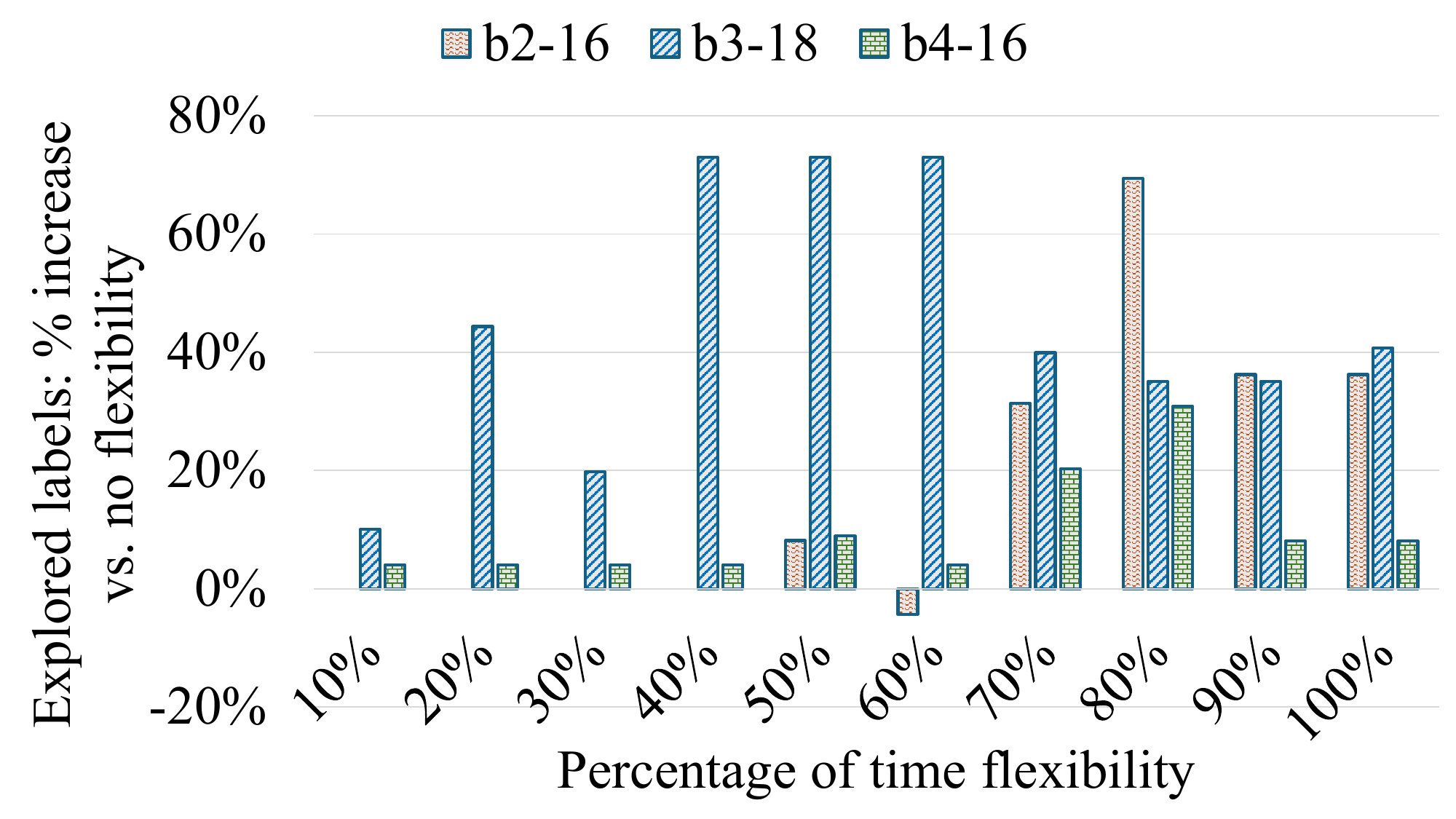}
    \caption{Growth of explored labels under increasing time flexibility}
    \label{fig:time_flexibility}
\end{figure}
To evaluate the price of time flexibility, we investigate how increasing the proportion of flexible time windows affects the number of explored labels during the pricing process. Figure~\ref{fig:time_flexibility} shows the percentage increase in label exploration for three benchmark instances, ``b2-16'', b3-18'', and ``b4-16'', as the time flexibility ratio increases from 10\% to 100\%. The baseline is the classical C-DARP model (0\% flexibility). The vertical axis measures the percentage increase in the number of explored labels compared to this baseline.

The results indicate that as more service nodes are permitted to deviate from their hard time windows, the solution space expands and leads to more labels being generated. However, our proposed dominance rules effectively mitigate this growth by pruning suboptimal or redundant paths, thus keeping the computational burden within manageable limits. From a system design perspective, this implies that moderate time flexibility can be introduced without incurring prohibitive computation costs.

\begin{figure}[htbp]
    \centering
    \begin{subfigure}[b]{0.48\textwidth}        \includegraphics[width=\linewidth]{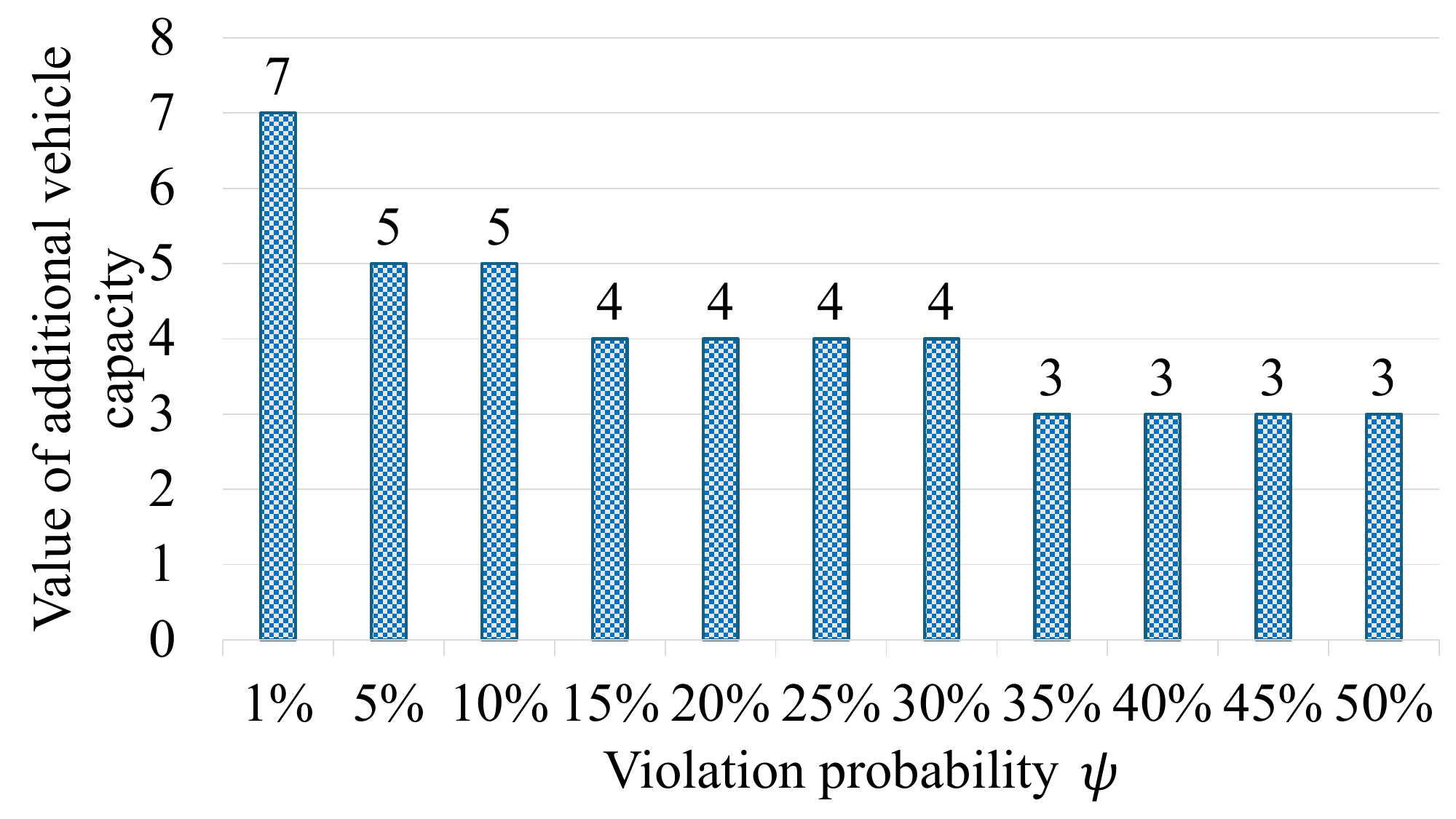}
        \caption{Increase in vehicle capacity required compared to baseline}
    \end{subfigure}
    \hfill
    \begin{subfigure}[b]{0.48\textwidth}
        \includegraphics[width=\linewidth]{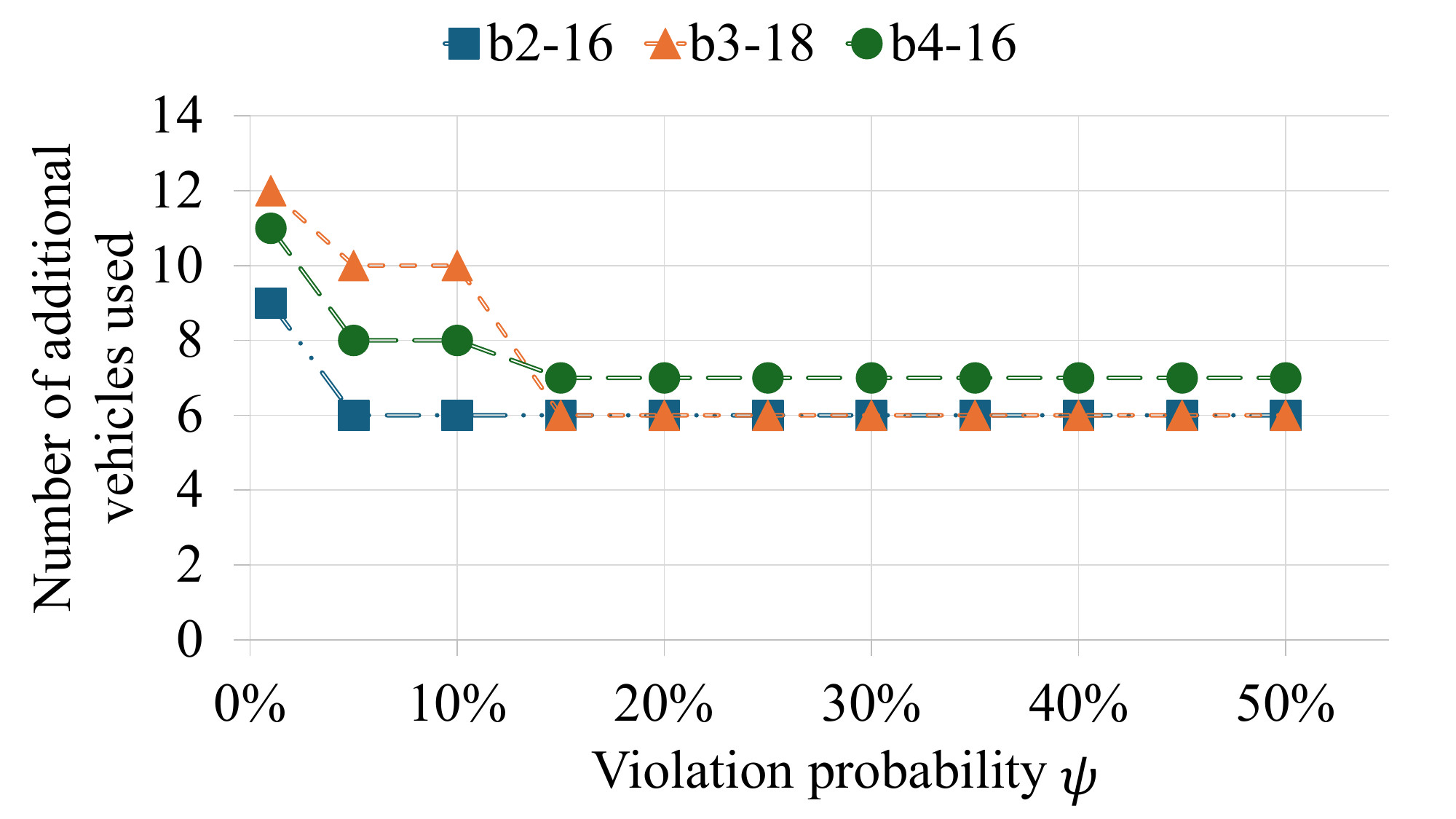}
        \caption{Increase in number of vehicles used compared to baseline}
    \end{subfigure}
    \caption{Resource costs of robustness under different violation probabilities}
    \label{fig:price_robustness}
\end{figure}

To evaluate the price of robustness, we examine the operational resources required to achieve different levels of service reliability, expressed through the violation probability $\psi$. Figure~\ref{fig:price_robustness} quantifies two key cost metrics under varying levels of $\psi$: the additional vehicle capacity required (Figure~\ref{fig:price_robustness}(a)) and the increase in fleet size (Figure~\ref{fig:price_robustness}(b)).

Figure~\ref{fig:price_robustness}(a) shows that reducing $\psi$ (i.e., demanding higher service reliability) significantly increase the vehicle capacity needed per vehicle. For example, maintaining a 95\% service level ($\psi=5\%$) requires increasing vehicle capacity from the baseline value 6 to 11, an 85\% increase. This reflects the need to preserve more empty seats to accommodate unexpected demand and reduce the risk of service failure.

Figure~\ref{fig:price_robustness}(b) reveals a similar trend in terms of the number of vehicles used. To meet a low violation probability, such as $\psi=5\%$, the system may require up to three times as many vehicles to compared to the baseline configuration. Although the resource requirements gradually decline as $\psi$ increases, the operational cost of robustness remains significant even at moderate reliability levels.

In summary, introducing time flexibility enlarges the solution space and increases label exploration, while robustness imposes substantial operational costs in terms of vehicle capacity and fleet size. Nevertheless, our BCP algorithm equipped with dominance rules and conservative planning effectively manages these trade-offs. These findings provide valuable insights for designing microtransit systems that balance customer flexibility, service reliability, and resource efficiency.

\subsection{Case study} \label{sec:case_study}
This section presents a case study based on transportation data\footnote{\url{https://data.cityofchicago.org/Transportation/Transportation-Network-Providers-Trips-2023-2024-/n26f-ihde/data_preview}} collected in Chicago between 2023 and 2024. We focus on the Austin community, a relatively underserved and geographically peripheral area, following studies in~\cite{ZUNIGAGARCIA2022100032,li2024flexible}. To simulate microtransit operations, we consider a two-hour planning horizon and evaluate three representative time periods: AM-PEAK (6:00-8:00), OFF-PEAK (10:00-12:00), and PM-PEAK (16:00-18:00). For each period, we optimize the proposed flexible and robust microtransit model using the modified BCP algorithm described in Section \ref{sec:branch_and_cut_and_price}. 

From the dataset, we extract relevant trip-level information, including pick-up and drop-off locations, pick-up and drop-off times, and the number of passengers, covering the one month period from December 1 to December 31, 2024. A parking facility located at coordinates $(41.8871,-87.7744)$ near the Austin subway station is designated as both the origin and destination depot. We set the following: travel distances and durations between any two location are computed using an open-source routing machine~\citep{luxen-vetter-2011}. Each vehicle is assumed to have a capacity of 10 passengers \citep{ghimire2024policy}, and the service time per request is assumed to be 30 seconds. To model demand uncertainty, we apply the sample average approximation method~\citep{patil2011sample} to the historical data to estimate the mean, variance, and bounded support (upper and lower bounds) of passenger load for each request, and assume the service level is $95\%$ (violation probability $\psi$=5\%). 

\subsection{Services comparison}
We conduct a comparative analysis of five service strategies to demonstrate the benefits of integrating time flexibility and robustness into microtransit systems. By comparing these strategies, we aim to highlight how time flexibility can improve service efficiency and reduce operational cost, while robustness enhances reliability under uncertain demand. The five strategies considered are: ride-hailing (RH''), traditional microtransit (M''), flexible microtransit (FM''), robust microtransit (RM''), and flexible and robust microtransit (``FRM''). Each strategy is evaluated across three representative time periods: AM-PEAK (6:00–8:00), OFF-PEAK (10:00–12:00), and PM-PEAK (16:00–18:00). Ride-hailing assumes proactive repositioning with full prior knowledge of all request locations and times, similar to the pre-booked nature of microtransit requests. To ensure consistency in comparison, travel time and travel mileage for microtransit operations exclude the initial travel from the depot to the first request, aligning with the assumptions in ride-hailing operations. All microtransit variants are solved using our modified BCP framework.

Table~\ref{tab:service_comparison} summarizes the key performance metrics. Column $N_r$'' denotes the number of requests served ; $V_u$'' is the number of vehicles used; ``TM'' indicates total travel mileage; ``T'' is total travel time (in minutes); ``CPU'' is computation time in seconds; ``WT'' is the average waiting time per request; ``TT'' is the total time (travel + total waiting time); and SR (\%)'' denotes service rate, which is computed as: 
$$
\text{SR}\%=\dfrac{\text{Total actually served requests}}{\text{Total requests}}
$$
The numerator counts the number of requests successfully served, i.e., a request is served successfully when all passengers of the request are served by the vehicle, otherwise, the vehicle will lose this request. The denominator is the number of total pre-booked requests. We simulate 50 scenarios where each request's actual load equals the pre-booked load plus a random offset $g \in \{-1,0,1,2,3,4\}$, each with equal probability. 
\begin{table}[ht!p]
    \centering
    \begin{tabular}{c|ccccccccc}
    \toprule
         & & $N_r$ & $V_u$ & TM & T & CPU & WT & TT & SR$\%$\\
    \hline
         \multirow{5}{*}{AM-PEAK}& RH & 21 & 20 & 151.61 & 233.29 & \textbf{0.0009} & \textbf{0.01} & 233.40 & 52.00\\
         & M & 21 & \textbf{7} & 124.62 & 190.38 & 0.0170 & 0.32 & 197.17 & 94.19 \\
         & FM & 21 & \textbf{7} & \textbf{112.33} & \textbf{172.57} & 0.030 & 0.51 & \textbf{183.24} & 94.38 \\
         & RM & 21 & 10 & 138.99 & 213.07 & 0.0150 & 0.09 & 214.96 & \textbf{96.52} \\ 
         & FRM & 21 &9 &124.88 & 185.56 & 0.0370 & 0.34 & 192.68 & \textbf{96.52} \\
    \hline
     \multirow{5}{*}{OFF-PEAK}& RH & 29 & 28 & 207.20 & 335.54 & \textbf{0.0007} & \textbf{0.01} & 335.78 & 45.20 \\
         & M & 29 & \textbf{12} & 169.36 & 251,72 & 0.0990 & 0.44 & 264.58 & 89.79\\
         & FM & 29 & \textbf{12} & \textbf{157.91} & \textbf{238.90} & 1.8830 & 0.72 & \textbf{259.86} & 91.90\\
         & RM & 29 & 15 & 180.58 & 267.55 & 0.0720  & 0.50 & 281.95 & 91.59 \\
         & FRM & 29 & 15 & 152.59 & 241.40 & 0.1030 & 0.60 & 258.72 & \textbf{93.66} \\
    \hline
     \multirow{5}{*}{PM-PEAK}& RH & 36 & 35 & 253.78 & 390.56 & \textbf{0.0008} & \textbf{0.01} &  390.92 & 52.61 \\
         & M & 36 & 14 & 168.53 & 264.21 & 0.1860 & 0.38 & 277.86 & 86.58\\
         & FM & 36 & \textbf{13} & \textbf{158.88} & \textbf{240.18} & 0.2490 & 0.60 & \textbf{261.85} & 87.78 \\
         & RM & 36 & 20 & 204.59 & 312.37 & 0.0550 & 0.61 & 334.18 & \textbf{96.42} \\
         & FRM & 36 & 17 & 188.64 & 259.31 & 0.1290 & 0.61 & 281.25 & \textbf{96.42}\\
    \bottomrule
    \end{tabular}
    \caption{Service comparison between different modes}
    \label{tab:service_comparison}
\end{table}
We draw the following key observations from Table~\ref{tab:service_comparison}:
\begin{enumerate}
    \item \textbf{Vehicle usage efficiency}: The number of vehicles used follows the trend: $\text{RH} \gg \text{RM} \geq \text{FRM} \geq \text{M} \geq \text{FM}$. Ride-hailing required the most vehicles due to its one-to-one request fulfillment model. Flexible microtransit (FM) achieves a 67.67\% reduction in vehicle usage compared to ride-hailing on average, by exploiting time flexibility to bundle more requests per vehicle. In contrast, robust microtransit (RM) demands 36.91\% more vehicles than traditional microtransit (M), as additional seat capacity is reserved to ensure reliability.
    \item \textbf{Time and mileage savings}: FM consistently achieves the lowest travel mileage (TM), travel time (T), and total time (TT) across all periods, demonstrating the efficiency benefits of time flexibility. On average, it saves 11.55 minutes in total time and 11.13 miles in travel distance compared to M. This highlights its value in reducing both passenger consumption time and operational cost.
    \item \textbf{Trade-off for robustness}: FRM requires a 26\% increase in vehicles compared to FM but offers significantly improved reliability. It achieves the highest average service rate (SR\%) of 95.53\%, effectively balancing flexibility and protection against demand uncertainty.
    \item \textbf{Service Reliability}: RM and FRM models substantially outperform others in SR\%, highlighting the importance of robustness for real-world operations. RH consistently underperforms in this metric, indicating that proactive rebalancing without robustness consideration is insufficient for uncertain demand.
\end{enumerate}
Microtransit services provide a more resource-efficient and reliable alternative to ride-hailing, particularly in environments with limited vehicle supply. Flexible microtransit enhances efficiency through time flexibility, while robust microtransit improves reliability at the cost of additional vehicle resources. The combined flexible and robust approach (FRM) offers the most balanced solution, achieving high service rates under uncertainty while maintaining operational feasibility.

\subsection{Sensitivity analysis}
To further evaluate the robustness and adaptability of our proposed microtransit model, we conduct a series of sensitivity analyses on three key dimensions: time flexibility, violation probability $\psi$, and demand uncertainty scenarios. While previous comparisons demonstrated the overall benefits of integrating flexibility and robustness, these additional experiments aim to quantify how varying levels of flexibility and conservativeness affect system performance and resource requirements. Specifically, we assess how different percentages of flexible time windows influence travel efficiency, how varying $\psi$ values impact the trade-off between service reliability and resource usage, and how alternative demand modeling assumptions shape the quality and feasibility of resulting schedules.

\paragraph{Sensitivity analysis on time flexibility}
We evaluate the impact of varying time flexibility levels on microtransit performance across three distinct time periods: AM-PEAK, OFF-PEAK, and PM-PEAK. The percentage of time flexibility ranges from 0\% to 100\%, where 0\% represents a fully rigid system in which all requests have hard (critical) time windows, and 100\% represents a fully flexible system where all requests allow relaxed scheduling through soft time windows. Intermediate levels, such as 30\% or 60\%, correspond to hybrid settings where a specified proportion of requests are randomly assigned as flexible. For each level, we conduct a sensitivity analysis based on 50 randomized instances reflecting the given percentage of flexible requests.
\begin{figure}[htbp]
    \centering
    \begin{subfigure}[b]{0.32\textwidth}        \includegraphics[width=\linewidth]{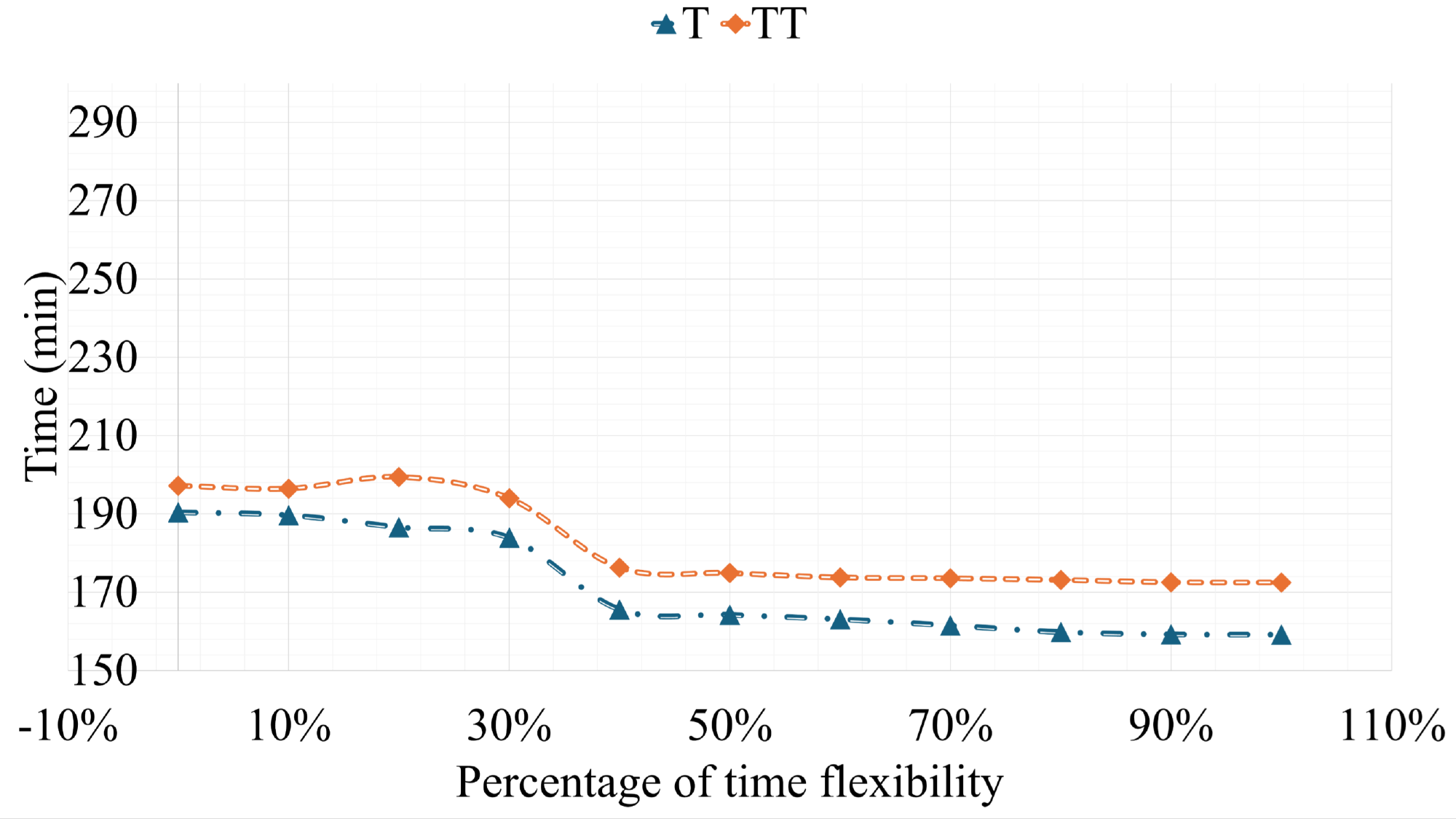}
        \caption{AM-PEAK}
    \end{subfigure}
    \hfill
    \begin{subfigure}[b]{0.32\textwidth}
        \includegraphics[width=\linewidth]{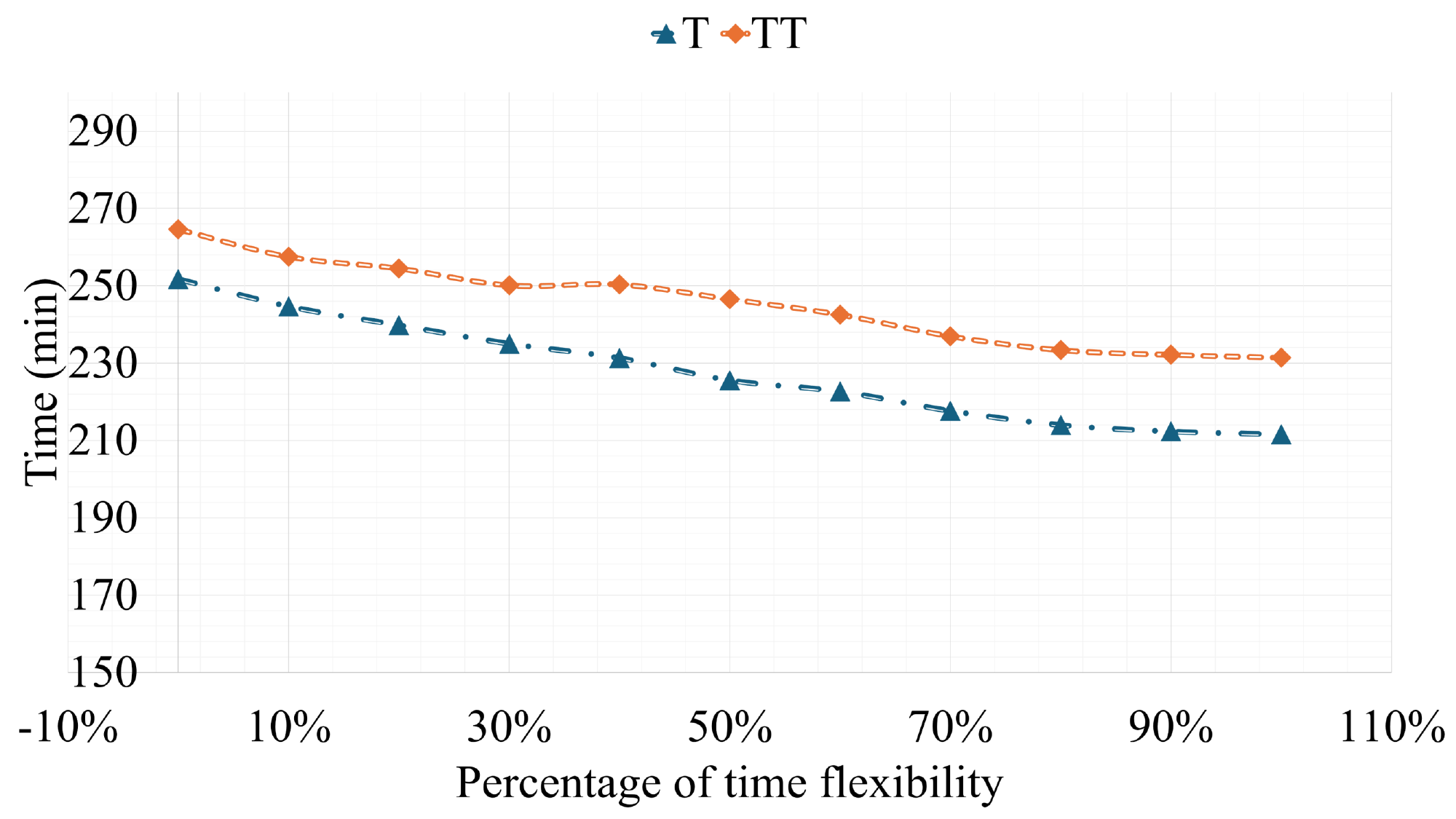}
        \caption{OFF-PEAK}
    \end{subfigure}
    \hfill
    \begin{subfigure}[b]{0.32\textwidth}
        \includegraphics[width=\linewidth]{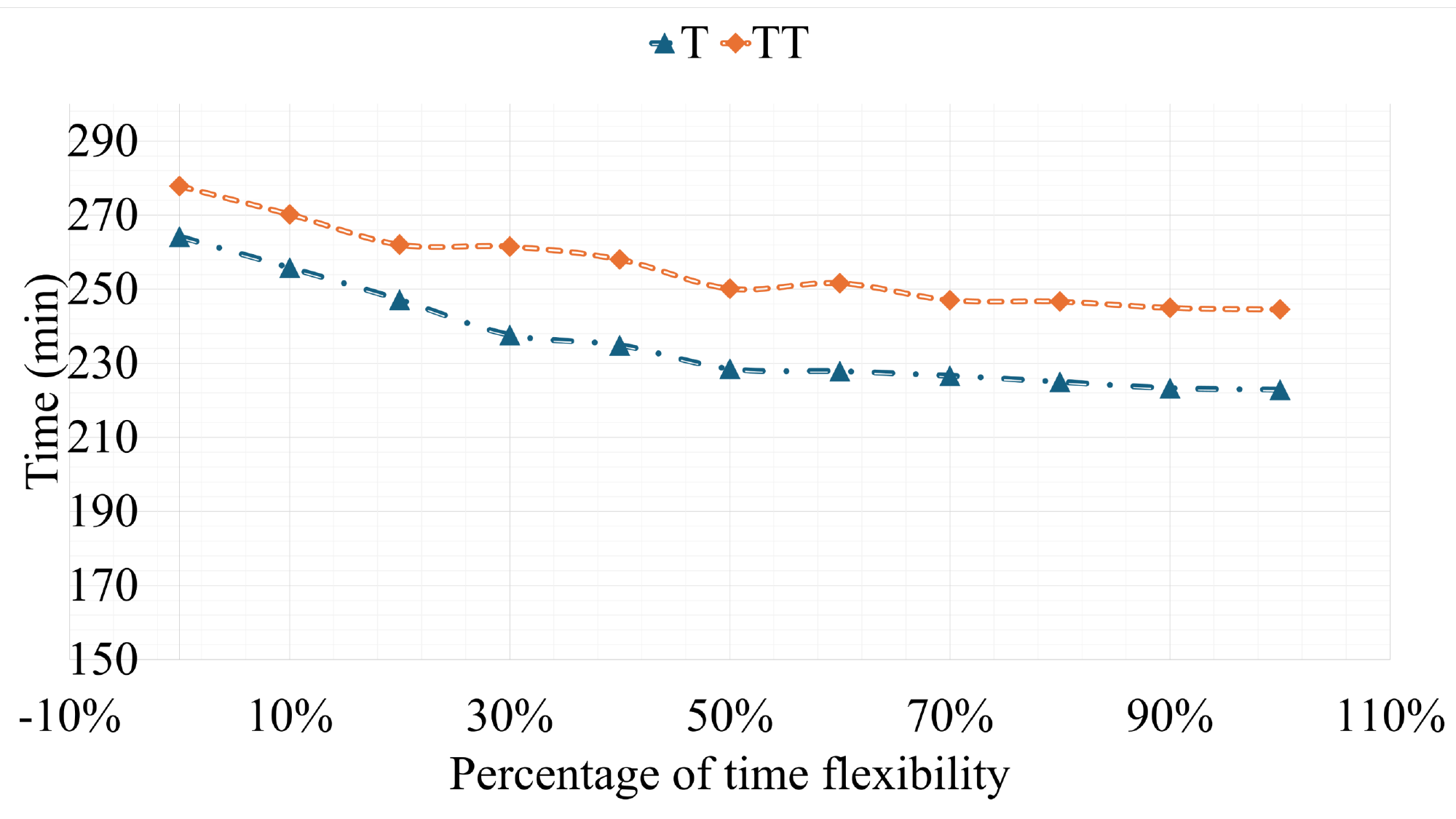}
        \caption{PM-PEAK}
    \end{subfigure}
    \caption{Sensitivity analysis on time flexibility percentage}
    \label{fig:time_flexibility}
\end{figure}

Figure~\ref{fig:time_flexibility} shows how total travel time (T) and total time (TT) evolve with increasing time flexibility. Across all three periods, we observe a clear downward trend: as more requests become time-flexible, both T and TT decrease significantly. The greatest improvements occur within the first 50\% increase in flexibility, after which the benefits continue at a slower, steadier rate.

These findings highlight the dual value of time flexibility. First, it enables the system to generate more efficient schedules with lower travel time, thereby reducing overall service costs. Second, it allows for greater customization in service delivery by accommodating varying degrees of passenger time sensitivity. From the platform's perspective, encouraging passengers to accept more flexible pick-up or drop-off windows can substantially reduce travel distances, minimize vehicle idle time, and enhance fleet utilization—ultimately translating into lower operational burden and improved service quality.

\paragraph{Sensitivity analysis on violation probability $\psi$} We conduct a sensitivity analysis on the violation probability parameter $\psi$, evaluating a range of values: $\psi \in \{1\%, 5\%, 10\%, 15\%, 20\%,25\%, 30\%, 35\%, 40\%, 45\%, 50\%\}$ across three time periods: AM-PEAK, OFF-PEAK, and PM-PEAK. For each setting, we simulate 50 randomized instances and evaluate two key performance metrics: the number of vehicles used ($N_v$) and the percentage of successfully served requests (SR\%). 
\begin{figure}[htbp]
    \centering
    \begin{subfigure}[b]{0.32\textwidth}
        \includegraphics[width=\linewidth]{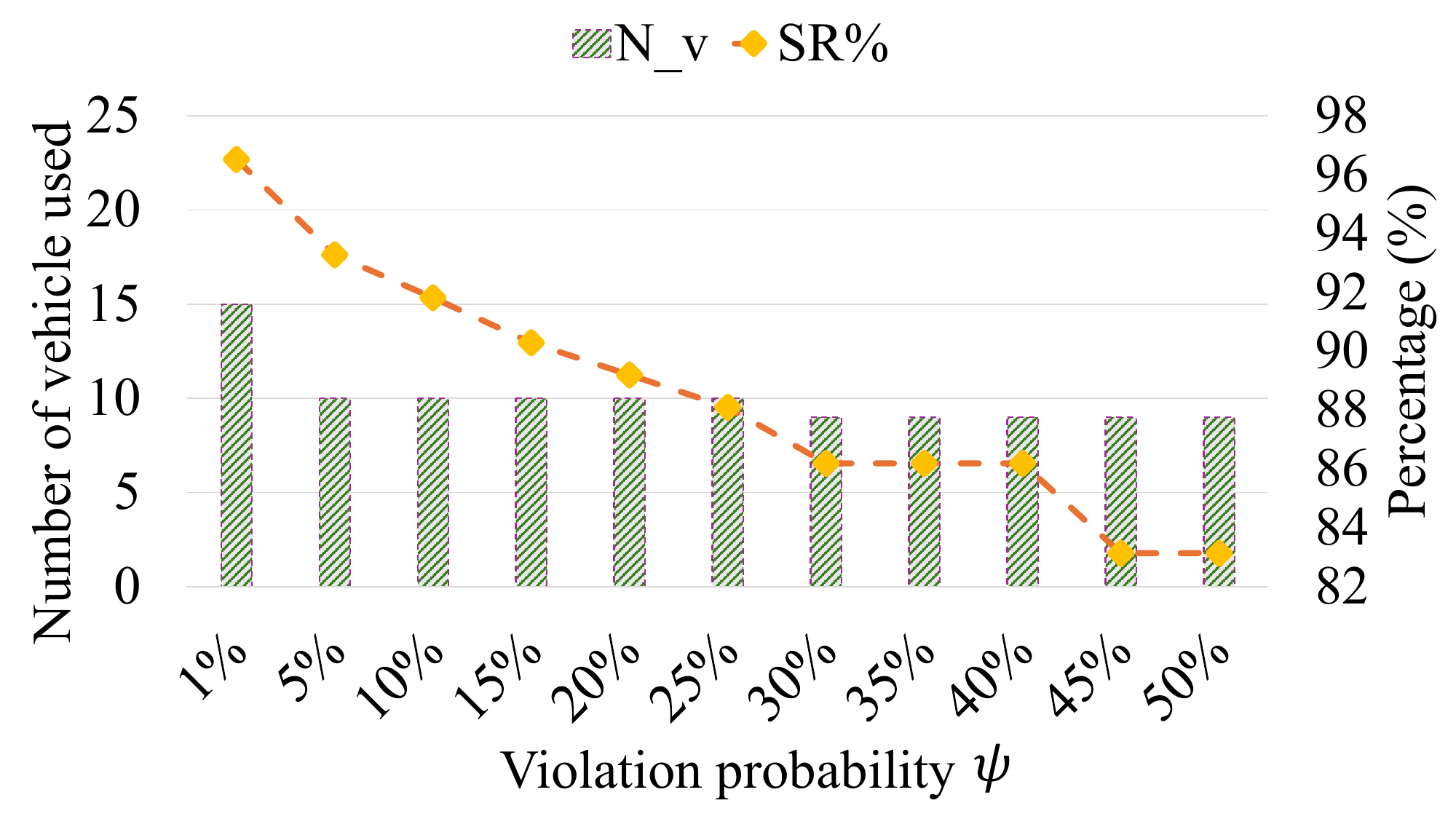}
        \caption{AM-PEAK}
    \end{subfigure}
    \hfill
    \begin{subfigure}[b]{0.32\textwidth}
        \includegraphics[width=\linewidth]{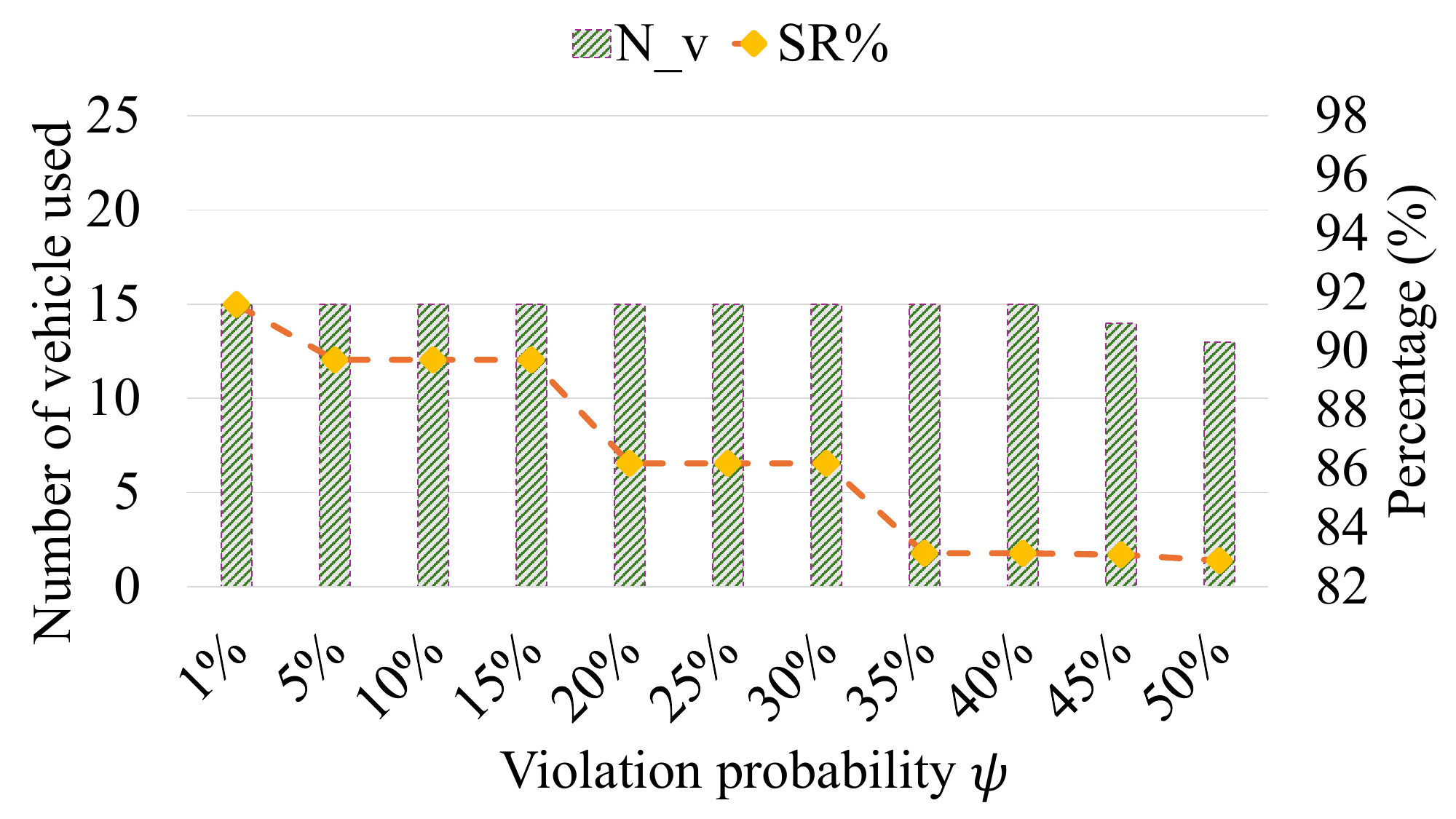}
        \caption{OFF-PEAK}
    \end{subfigure}
    \hfill
    \begin{subfigure}[b]{0.32\textwidth}
        \includegraphics[width=\linewidth]{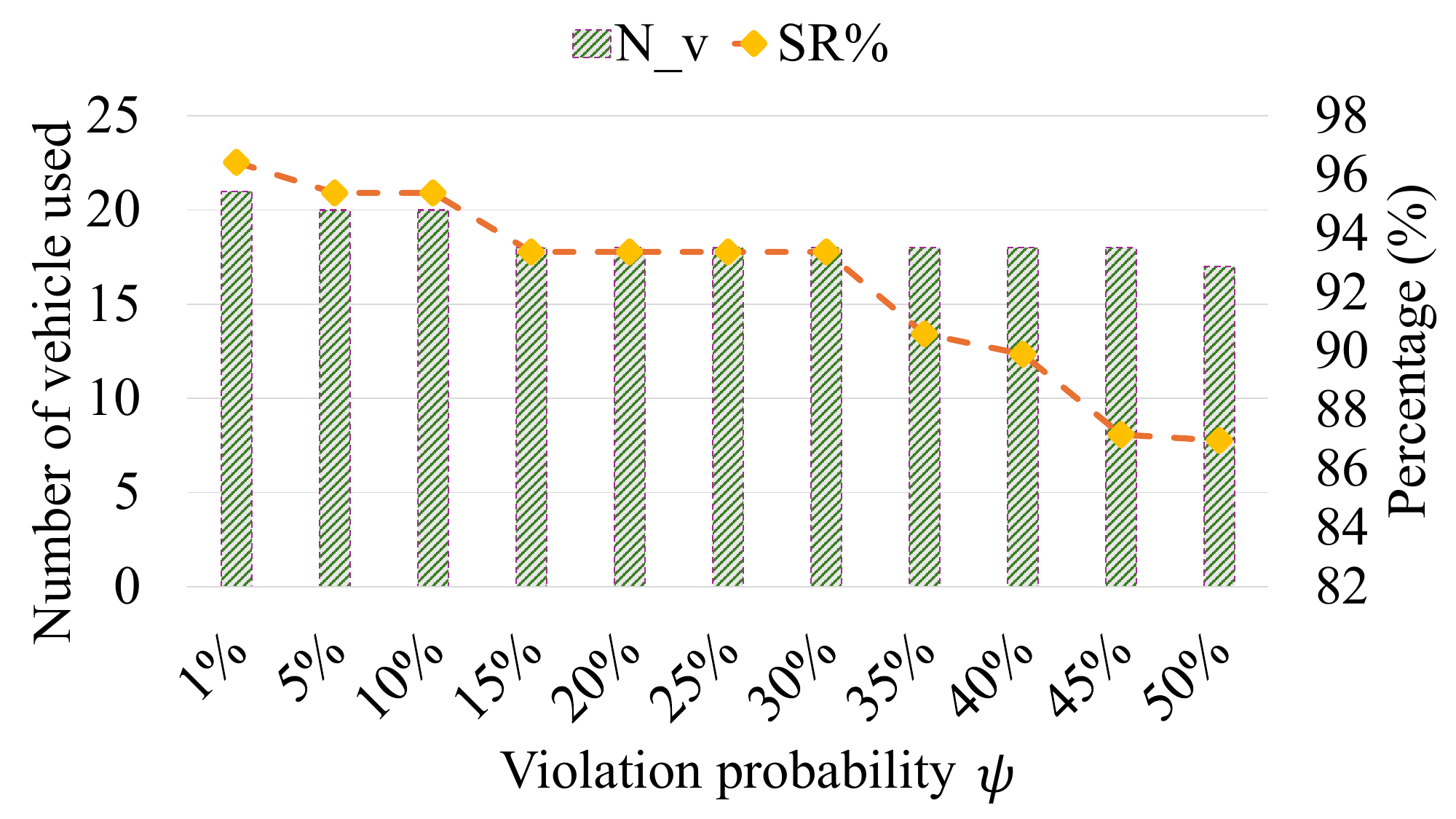}
        \caption{PM-PEAK}
    \end{subfigure}
    \caption{Sensitivity analysis on violation probability $\psi$}
    \label{fig:violation probability}
\end{figure}
Figure~\ref{fig:violation probability} illustrates the trade-off between service reliability and resource usage as $\psi$ varies. As expected, increasing the violation probability allows the system to operate under less conservative assumptions. We observe that the service rate (SR\%) decrease sharply as $\psi$ increases (shown in AM-PEAK and PM-PEAK). In contrast, the number of vehicles used decrease only slightly, indicating that higher violation probabilities do not substantially increase resource requirements.

These results suggest that raising $\psi$ can lead to better vehicle utilization without significantly compromising operational costs, but at the expense of service reliability. In practical terms, selecting a moderately high $\psi$ (e.g., between 5\% and 10\%) may strike a desirable balance between maintaining service quality and limiting fleet size. This flexibility allows operators to turn system conservativeness based on available resources and acceptable levels of risk.

\paragraph{Sensitivity analysis on demand scenarios} We further compare five demand modeling scenarios to evaluate their impact on service quality and reliability: joint independent Normal distribution (JIND), joint dependent Normal distribution (JDND), ambiguity set (AS) based on distributionally robust optimization with known mean and variance, bounded support (BS) as described in Section~\ref{sec:bounded_support}, and a deterministic case (D). For detailed formulations of the JIND, JDND, and AS models, we refer readers to \cite{dinh2018exact}.

\begin{table}[h]
    \centering
    \begin{tabular}{c|ccccccccc}
    \toprule
         & & $N_r$ & $V_u$ & TM & T & CPU (s) & WT & TT & SR\% \\
    \hline
         \multirow{5}{*}{AM-PEAK} & JIND & 21 & 8 & 186.18 & 277.96 & \textbf{0.026} & 0.15 & \textbf{194.38} & 88.86 \\
         & JDND & 21 & 13 & 232.51 & 344.02 & 444.399 & 0.08 & 198.22 & 93.89\\
         & AS & 21 & 21 & 428.45 & 405.78 & 27.076 & \textbf{0.01} & 232.05 & 94.26 \\
         & BS & 21 & 10 & 219.40 & 266.2 & 0.054 & 0.09 & 214.96 & \textbf{96.19} \\
         & D & 21 & \textbf{7} & \textbf{172.25} & \textbf{260} & 0.03 & 0.32 & 197.10 & 84.96 \\
    \hline
     \multirow{5}{*}{OFF-PEAK} & JIND & 29 & \textbf{12} & 292.22 & 417.3 & \textbf{0.13} & 0.54 & 267.75 & 89.62\\
         & JDND & 29 & 17 & 352.00 & 497.72 & 182.97 & 0.24 & 296.15 & 94.38 \\
         & AS & 29 & 29 & 427.71 & 595.18 & 47.735 & \textbf{0.01} & 337.27 & 96.14 \\
         & BS & 29 & 15 & 327.98 & 466.69 & 0.105 & 0.5 & 282.05 & \textbf{96.59} \\
         & D & 29 & \textbf{12} & \textbf{286.26} & \textbf{410.78} & 0.03 & 0.44 & \textbf{264.48} & 89.62\\
    \hline
          \multirow{5}{*}{PM-PEAK} & JIND & 36 & 15 & 294.98 & 443.45 & 1.081 & 0.35 & 278.58 & 89.22\\
         & JDND & 36 & 22 & 372.69 & 553.27 & 465.69 & 0.14 & 353.95 &94.89  \\
         & AS & 36 & 36 & 473.29 & 687.54 & 86.247 & \textbf{0.01} & 392.83 & 95.89 \\
         & BS & 36 & 20 & 375.71 & 555.53 & 0.056 & 0.61 & 334.33 & \textbf{96.61} \\
         & D & 36 & \textbf{13} & \textbf{271.89} & 418.68 & \textbf{0.249} & 0.38 & \textbf{277.89} & 86.81\\
    \bottomrule
    \end{tabular}
    \caption{Robustness under different scenarios}
    \label{tab:robustness}
\end{table}
Table~\ref{tab:robustness} presents the comparative results using the following performance indicators: number of served requests ($N_r$), number of vehicles used ($N_u$), total mileage (TM), total travel time (T), average waiting time (WT), total time (TT = travel + waiting), CPU time, and service rate (SR\%).

Our findings reveal several important insights. The bounded support (BS) scenario achieves the highest service reliability, as reflected by the highest SR\% of 96.46\% on average, while also maintaining a reasonable number of vehicles used and moderate CPU time. This suggests that BS offers a balanced and practical approach to modeling robustness under demand uncertainty.

The deterministic (D) case produces the lowest SR\%, as it does not account for demand variability. It achieves the shortest travel time and mileage and uses the fewest vehicles, but sacrifices reliability significantly. Interestingly, the JIND case produces almost identical schedules to the deterministic case, indicating that assuming independent Gaussian demand fails to capture the variability needed to support robust service delivery.

In contrast, the JDND scenario improves SR\% by modeling correlations among requests, resulting in better reliability. However, this comes at the cost of higher vehicle usage, longer CPU time, and increased travel metrics. The ambiguity set (AS) case is the most conservative, often assigning one vehicle per request, which results in the highest operational cost and extremely low system efficiency despite achieving high SR\%.

In summary, among the stochastic models considered, the bounded support (BS) scenario is the most ideal for capturing demand uncertainty. It achieves the best trade-off between reliability, resource usage, and computational efficiency, making it a practical and robust choice for real-world microtransit planning.

\section{Conclusion}\label{sec:conclusion}
This paper presents a comprehensive framework for designing a flexible and robust microtransit system that jointly addresses two fundamental operational challenges: passenger time sensitivity and demand uncertainty. We formalize this problem as the Chance-Constrained Dial-a-Ride Problem with Soft Time Windows (CCDARP-STW), which integrates soft time windows to model heterogeneous urgency among passengers and employs chance constraints to manage demand uncertainty. To enhance tractability, we reformulate the problem into a robust trip-based model and apply a linear relaxation using Hoeffding's inequality, requiring only bounded-support statistics.

To solve CCDARPSTW efficiently, we develop a customized Branch-and-Cut-and-Price (BCP) algorithm that incorporates a forward labeling scheme supporting soft time windows and a novel probabilistic dominance rule (PDR) for label pruning under uncertainty. Computational experiments on benchmark instances confirm that the proposed algorithm maintains scalability even as complexity increases with time flexibility and robustness features. Notably, the PDR yields a 17.40\% reduction in label exploration and a 22.27\% reduction in CPU time, illustrating the value of risk-aware pruning in stochastic routing. Beyond this specific context, the PDR is generalizable to a wide range of robust combinatorial optimization problems that involve probabilistic feasibility constraints. 

A case study using real-world ride-hailing data from the Austin neighborhood in Chicago further demonstrates the practical value of the proposed approach. Flexible microtransit achieves up to 67.67\% reduction in vehicle usage and notable time and mileage savings compared to conventional microtransit. When robustness is introduced to guarantee a 5\% violation probability, vehicle usage increases modestly by 26\%, but service reliability reaches 95.53\% on average, with the bounded-support method outperforming normal distribution and ambiguity-set-based models in both effectiveness and computational efficiency. 

Finally, sensitivity analyses confirm that (1) moderate time flexibility delivers substantial efficiency benefits; (2) moderate violation probabilities provide a desirable balance between resource use and reliability; and (3) the bounded-support approach offers the best trade-off among stochastic demand models, delivering high service rates (up to 96.61\%) with limited computational and operational usages. In sum, this paper provides a unified and scalable optimization framework for real-world microtransit scheduling under demand uncertainty, enabling transportation agencies to adapt services to diverse user needs while ensuring reliability. 

Three promising directions for future work remain. First, the current model focuses on unexpected requests that share the same pickup and drop-off locations as pre-booked trips. Future research could extend this by considering more diverse forms of unexpected demand, including requests with new origins or destinations and varying time windows. Second, developing a unified scheduling platform that simultaneously manages both real-time and pre-booked requests could further improve system responsiveness and resource efficiency. Finally, while this study models demand uncertainty using a bounded-support approach, future work could explore more realistic uncertainty structures by embedding demand forecast errors or systematic deviations, better reflecting the operational challenges faced in practice.

\section*{Acknowledgement}
This research is funded by the U.S. Department of Energy, Office of Energy Efficiency for the Co-E3T (Energy-Efficient and Equitable Transit through user-centric hardware and software Co-development and community Co-design) project.

\bibliographystyle{plainnat}
\bibliography{sample}

\appendix
\section{Notations} \label{sec:appendix_notation}

\begin{table}[ht!p]
    \centering
    \begin{tabular}{p{2cm}|p{13cm}}
    \toprule
    \multicolumn{2}{c}{Parameters Notation}\\
    \hline
         $n$ & the number of requests;\\
         $\mathcal{P}$ & the set of pick-up nodes and $\mathcal{P}=\{1,2,\cdots,n\}$; \\
         $\mathcal{D}$ & the set of drop-off nodes and $\mathcal{D}=\{n+1,n+2,\cdots,2n\}$; \\
         $P$ & the set of flexible pick-up nodes and $P \subseteq \mathcal{P}$, therefore $\mathcal{P}\setminus P$ is the set of time-sensitive pick-up nodes; \\
         $D$ & the set of flexible drop-off nodes and $D \subseteq \mathcal{D}$, therefore $\mathcal{D}\setminus D$ is the set of time-sensitive pick-up nodes;\\
         $\mathcal{G}=(\mathcal{N},\mathcal{A})$ & the transportation network where $\mathcal{N}=\mathcal{P}\cup \mathcal{D} \cup \{0,2n+1\}$ and $0$,$2n+1$ are origin depot and destination depot, respectively; \\
         $[e_i,l_i]$ & time window for nodes $i \in \mathcal{N}$, $e_i$ represents the earliest service time and $l_i$ represents the latest service time;\\
         $1-\psi$ & the confidence level derived from the demand chance constraint; \\
          $\mathcal{K}$ & the set of the vehicle $k \in \mathcal{K}$;\\
          $M$ & the capacity of the vehicle $k \in \mathcal{K}$;\\
          $B_i$ & the maximum duration time for each request $i$ \\
          $T$ & the maximum time horizon \\
          $c_{ij}$ & the travel cost for arc $(i,j) \in \mathcal{A}$;\\
          $H_i^k(\cdot)$ & a generalized penalty cost function of vehicle $k$ serving node $i$.\\
          \hline
          \multicolumn{2}{c}{Random Variables}\\
          \hline
          $\tilde{w}_i(\omega)$ & the random demand at node $i \in \mathcal{P}$ based on scenarios $\omega \in \Omega$ and $\tilde{w}_i(\omega)=-\tilde{w}_{i+n}(\omega)$\\
          \hline
          \multicolumn{2}{c}{Decision Variables}\\
          \hline
          $x_{ij}^k$ & a binary variable, it is 1 if vehicle $k$ travels the arc $(i,j) \in \mathcal{A}$;\\
          $y_i^k$ & a continuous variable for the delay time at the node $i$ for vehicle $k$;\\
          $L_i^k$ & a continuous variable recording the number of onboarding customers of vehicle $k$ after visiting node $i$; \\
          $T_i^k$ & a continuous variable representing the visiting time at node $i$ of vehicle $k$.\\
          $D_i^k$ & a continuous variable of ride time for request $i$ on the vehicle $k$ \\
    \bottomrule
    \end{tabular}
    \caption{Notations}
    \label{tab:my_notation}
\end{table}

\section{Proof of Proposition \ref{proposition:demand_linear}}\label{appendix:linear_constraint}
Let $\tilde{w}_j$ denote the uncertain demand for request $j$, with known expectation $\mu_j=\mathbb{E}[\tilde{w}_j]$, and bounded support $a_j \leq \tilde{w}_j \leq b_j$. Suppose the requests in the onboard set $O_i$ are either independent or weakly dependent. According to Hoeffding's inequality (or its extension to weakly dependent variables as shown in \cite{janson2004large}), we have
    $$
    \text{Prob}\Bigg\{\sum_{j \in O_i} \tilde{w}_j -\sum_{j \in O_i}\mu_j \geq t \Bigg \} \leq \text{exp}\Bigg(-\dfrac{2t^2}{\sum_{j \in O_i}(b_j-a_j)}\Bigg).
    $$
    To ensure the chance constraint
    $$
    \text{Prob}\Bigg\{\sum_{j \in O_i} \tilde{w}_j \leq M\Bigg\} \geq 1-\psi
    $$
    we need: 
    $$
    \text{Prob}\Bigg\{\sum_{j \in O_i}\tilde{w}_j-\sum_{j \in O_i} \mu_j \leq M - \sum_{j \in O_i} \mu_j\Bigg\} \leq \psi. 
    $$
    Let $t:=\sum_{j \in O_i}\mu_j -M$, then by Hoeffding's inequality, we must have:
    $$
    \psi \geq \text{exp}\Bigg(-\dfrac{2(\sum_{j \in O_i}\mu_j -M)^2}{\sum_{j \in O_i}(b_j-a_j)^2}\Bigg)
    $$
    Take the logarithm on both sides:
    $$
    -\log (\psi) \leq \dfrac{2(\sum_{j \in O_i}\mu_j-M)^2}{\sum_{j \in O_i})(b_j-a_j)^2}. 
    $$
    Rewriting: 
    $$
    \sum_{j \in O_i} \mu_j -M \leq \sqrt{\dfrac{1}{2}(\sum_{j \in O_i}(b_j-a_j)^2)\cdot \log(\dfrac{1}{\psi})}
    $$
   Now, to obtain a linear (safe) approximation, we apply the norm inequality:
   $$
   \sqrt{\sum_{j \in O_i}(b_j-a_j)^2} \leq \sum_{j \in O_i}(b_j-a_j), 
   $$
   which yields:
    $$
    \sum_{j \in O_i} \mu_j -M \leq  \sqrt{\log (1/\psi)/2}\cdot \sum_{j \in O_i} (b_j-a_j).
    $$
   Rearanging terms:
    $$
    \sum_{j \in O_i}\Big(\mu_j+\sqrt{\log (1/\psi)/2}\cdot(b_j-a_j)\Big) \leq M.
    $$
    This completes the proof of the conservative linear approximation. 

\section{Proof of Proposition \ref{proposition:strong_dominance}}\label{appendix:proof_dominance}
The proposed strong dominance rule is designed to prune the label space efficiently while accounting for cost, feasibility, and robustness under demand uncertainty.

\textbf{Condition 1:} $(z_1 \leq z_2, E_1 \leq E_2, V_1 \subseteq V_2, O_1 \subseteq O_2)$ captures the standard deterministic dominance criteria commonly used in labeling algorithms for shortest path problems with resource constraints~\citep{gschwind2018bidirectional}. A label with lower reduced cost and earlier arrival time offers more temporal slack and cost efficiency. Furthermore, if label $u_1$ has completed and opened a subset of the requests handled by $u_2$, it retains more routing flexibility and future extension potential.

\textbf{Condition 2}: This condition ensures that the feasibility of ride time constraints is preserved, as discussed in~\cite{Gschwind_2015}. The term $R^0(E)-E$ represents the time spent from pickup to drop-off (i.e., ride duration) for a request. A larger value of this expression for $u_1$ compared to $u_2$ implies that label $u_1$ allows a longer feasible ride time margin than $u_2$ for each active request. Similarly, the condition $R^0_1(B_1^0)\geq R_2^0(B_2^0)$ ensures that the drop-off time remains within acceptable bounds across all open requests, even under worst-case (latest) departure times. These comparisons guarantee that the ride time requirements, critical for service quality and feasibility, are at least as well satisfied in $u_1$ as in $u_2$.

\textbf{Condition 3:} introduces robustness into the dominance rule by comparing the violation probabilities $\Gamma_1$ and $\Gamma_2$, derived using Hoeffding's inequality. A smaller $\Gamma$ implies a lower chance of violating the capacity constraint under uncertain demand. By requiring $\Gamma_1 \leq \Gamma_2$, we prioritize solutions that are probabilistically more reliable.

This incorporation of a probabilistic dominance criterion reflects a deliberate modeling choice: to sacrifice a small amount of optimality in exchange for improved robustness. This is consistent with the philosophy of price of robustness~\citep{bertsimas2004price}, as discussed in robust optimization literature. In practice, a label with slightly higher reduced cost (and hence potentially worse objective value) may be retained if it leads to a lower risk of infeasibility. This trade-off acknowledges that solutions in real-world settings must remain feasible under uncertainty, even if that means accepting marginally suboptimal costs.

Thus, the rule preserves labels that dominate others not only in deterministic terms but also in their probabilistic resilience, thereby guiding the labeling algorithm to construct robust, feasible trips with acceptable performance loss.

\end{document}